# Some applications of the Stieltjes constants

Donal F. Connon

dconnon@btopenworld.com

14 January 2009

**Abstract**

In this paper we present some applications of the Stieltjes constants including, for example, new derivations of Binet's formulae for $\log \Gamma(x)$ and the evaluation of some integrals related to the Barnes multiple gamma functions. Several connections with the sine and cosine integrals are also depicted.

This paper is dedicated to the memory of the Irish mathematician, John Todd (1911-2007). Todd was one of the pioneers of numerical analysis who, inter alia, computed numerical approximations of the first 20 Stieltjes constants [36] in 1972.

## 1. Introduction

The generalised Euler-Mascheroni constants $\gamma_n$ (or Stieltjes constants) are the coefficients of the Laurent expansion of the Riemann zeta function $\varsigma(s)$ about $s = 1$

$$\text{(1.1)} \qquad \varsigma(s) = \frac{1}{s-1} + \sum_{n=0}^{\infty} \frac{(-1)^n}{n!} \gamma_n (s-1)^n$$

Since $\lim_{s \to 1} \left[ \varsigma(s) - \frac{1}{s-1} \right] = \gamma$ it is clear that $\gamma_0 = \gamma$. It may be shown, as in [32, p.4], that

$$\text{(1.2)} \qquad \gamma_n = \lim_{N \to \infty} \left[ \sum_{k=1}^{N} \frac{\log^n k}{k} - \frac{\log^{n+1} N}{n+1} \right] = \lim_{N \to \infty} \left[ \sum_{k=1}^{N} \frac{\log^n k}{k} - \int_1^N \frac{\log^n t}{t} dt \right]$$

where, throughout this paper, we define $\log^0 1 = 1$.

In 1985, using contour integration, Ainsworth and Howell [6] showed that for $n \geq 1$

$$\text{(1.3)} \qquad \gamma_n = 2 \operatorname{Re} \int_0^{\infty} \frac{(x-i) \log^n (1-ix)}{(1+x^2)(e^{2\pi x} - 1)} dx$$

and for $n = 0$ we have

$$\text{(1.4)} \qquad \gamma_0 = \gamma = \frac{1}{2} + 2 \operatorname{Re} \int_0^{\infty} \frac{x-i}{(1+x^2)(e^{2\pi x} - 1)} dx = \frac{1}{2} + 2 \int_0^{\infty} \frac{x}{(1+x^2)(e^{2\pi x} - 1)} dx$$

The reason for the additional factor of ½ in (1.4) becomes apparent upon examination of equation (1.5) below when we let $n = 0$ and $u = 1$. It appears that Ainsworth and Howell [6] indirectly employed the Abel-Plana summation formula (3.1) in their derivation of equation (1.3).

We may also note that for $n \geq 1$

$$\gamma_n = 2 \operatorname{Re} \int_0^{\infty} \frac{(x-i) \log^n (1-ix)}{(1+x^2)(e^{2\pi x} - 1)} dx = -2 \frac{d^n}{dt^n} \operatorname{Re} \int_0^{\infty} \frac{i(1-ix)^t}{e^{2\pi x} - 1} dx \bigg|_{t=-1}$$

More recently, Coffey [20] has extended the result in (1.3) to show that the Stieltjes constants may be represented by the integral

$$\text{(1.5)} \qquad \gamma_n(u) = \frac{1}{2u} \log^n u - \frac{1}{n+1} \log^{n+1} u - 2 \operatorname{Re} \int_0^{\infty} \frac{i(u+ix) \log^n (u-ix)}{(u^2 + x^2)(e^{2\pi x} - 1)} dx$$

where the Stieltjes constants $\gamma_n(u)$ are the coefficients in the Laurent expansion of the Hurwitz zeta function $\varsigma(s, u)$ about $s = 1$

$$\text{(1.6)} \qquad \varsigma(s,u) = \sum_{n=0}^{\infty} \frac{1}{(n+u)^s} = \frac{1}{s-1} + \sum_{n=0}^{\infty} \frac{(-1)^n}{n!} \gamma_n(u)(s-1)^n$$

and $\gamma_0(u) = -\psi(u)$, where $\psi(u)$ is the digamma function which is the logarithmic derivative of the gamma function $\psi(u) = \frac{d}{du} \log \Gamma(u)$. It is easily seen from the definition of the Hurwitz zeta function that $\varsigma(s,1) = \varsigma(s)$ and accordingly that $\gamma_n(1) = \gamma_n$. A different derivation of (1.5) is shown below in (4.2).

## 2. New derivations of Binet's formulae for $\log \Gamma(u)$

In 2007, it was shown in [22] that the Stieltjes constants could be represented by the logarithmic expansion

$$\text{(2.1)} \qquad \gamma_n(u) = -\frac{1}{n+1} \sum_{i=0}^{\infty} \frac{1}{i+1} \sum_{j=0}^{i} \binom{i}{j} (-1)^j \log^{n+1}(u+j)$$

A simplified proof of this expression is given below in (4.8). Using this formulation, a number of identities involving the Stieltjes constants were derived, including the following

$$\text{(2.2)} \qquad \frac{1}{2} \log(2\pi) - \log \Gamma(u) = 1 + \sum_{n=0}^{\infty} \frac{\gamma_{n+1}(u)}{n!}$$

In fact (2.2) may be obtained more directly by differentiating (1.6) with respect to $s$, letting $s=0$ and then deploying Lerch's identity (5.4).

We now substitute Coffey's representation (1.5) in (2.2) to obtain

(2.3)

$$\sum_{n=0}^{\infty} \frac{\gamma_{n+1}(u)}{n!} = \frac{1}{2u} \sum_{n=0}^{\infty} \frac{1}{n!} \log^{n+1} u - \sum_{n=0}^{\infty} \frac{1}{n!} \frac{1}{n+2} \log^{n+2} u - 2\,\text{Re} \int_0^{\infty} \frac{i(u+ix) \sum_{n=0}^{\infty} \frac{1}{n!} \log^{n+1}(u-ix)}{(u^2+x^2)(e^{2\pi x}-1)} dx$$

and, for the first summation, we easily see from the exponential series that

$$\sum_{n=0}^{\infty} \frac{1}{n!} \log^{n+1} u = u \log u$$

As regards the second summation we note that

$$\sum_{n=0}^{\infty} \frac{y^{n+2}}{(n+2)!} = e^y - 1 - y$$

and thus

$$\sum_{n=0}^{\infty} \frac{y^{n+1}}{(n+2)!} = \frac{e^y - 1 - y}{y}$$

Differentiation results in

$$\sum_{n=0}^{\infty} \frac{1}{n!} \frac{1}{n+2} y^n = \frac{d}{dy}\left[\frac{e^y - 1 - y}{y}\right] = \frac{e^y(y-1)+1}{y^2}$$

and we end up with

$$\sum_{n=0}^{\infty} \frac{1}{n!} \frac{1}{n+2} y^{n+2} = e^y(y-1)+1$$

which gives us

$$\sum_{n=0}^{\infty} \frac{1}{n!} \frac{1}{n+2} \log^{n+2} u = u(\log u - 1) + 1$$

Finally we have the exponential series

$$\sum_{n=0}^{\infty} \frac{1}{n!} \log^{n+1}(u - ix) = (u - ix)\log(u - ix)$$

and we then obtain for the real part of (2.3)

$$\sum_{n=0}^{\infty} \frac{\gamma_{n+1}(u)}{n!} = \frac{1}{2}\log u - u(\log u - 1) - 1 - 2\int_0^{\infty} \frac{\tan^{-1}(x/u)}{e^{2\pi x} - 1} dx$$

Then, using (2.2), we have determined Binet's second formula for $\log \Gamma(u)$ (which is derived in a very different manner in, for example, [49, p.251])

(2.4) $$\log \Gamma(u) = \left(u - \frac{1}{2}\right)\log u - u + \frac{1}{2}\log(2\pi) + 2\int_0^{\infty} \frac{\tan^{-1}(x/u)}{e^{2\pi x} - 1} dx$$

This formula was also derived by Ramanujan [10, Part II, p.221] in the case where $u$ is a positive integer.

The substitution $x = ut$ gives us

$$\log \Gamma(u) = \left(u - \frac{1}{2}\right) \log u - u + \frac{1}{2}\log(2\pi) + 2u \int_0^\infty \frac{\tan^{-1} t}{e^{2\pi ut} - 1} dt$$

We now use (2.4) to derive Binet's first formula for $\log \Gamma(u)$. For convenience we define $f(u, x)$ by

$$f(u, x) = \int_0^\infty \frac{e^{-uy} \sin(xy)}{y} dy$$

and differentiation results in

$$\frac{\partial}{\partial u} f(u, x) = -\int_0^\infty e^{-uy} \sin(xy) \, dy$$

(alternatively we could have considered the other partial derivative $\frac{\partial}{\partial x} f(u, x)$).

The substitution $t = xy$ gives us

$$\int_0^\infty e^{-uy} \sin(xy) \, dy = \frac{1}{x} \int_0^\infty e^{-\frac{u}{x} y} \sin t \, dt$$

and using the familiar Laplace transform

(2.5)  $$\int_0^\infty e^{-at} \sin t \, dt = \frac{1}{1 + a^2}$$

we see that

$$\frac{\partial}{\partial u} f(u, x) = -\frac{x}{u^2 + x^2}$$

Integration of this equation gives us

$$f(u, x) = -\tan^{-1}(u/x) + c$$

Letting $u \to \infty$ we see from the defining integral that $f(u, x) \to 0$ and thus $c = \pi/2$.

Since

$$\tan^{-1} a + \tan^{-1} b = \tan^{-1}\left(\frac{a + b}{1 - ab}\right)$$

we see that

(2.6) $$\tan^{-1}(x/u) = \frac{\pi}{2} - \tan^{-1}(u/x) \text{ for } u/x > 0$$

We then obtain the well-known integral

(2.7) $$\tan^{-1}(x/u) = \int_0^\infty \frac{e^{-uy} \sin(xy)}{y} dy$$

Rigorous derivations of (2.7) are contained in [8, p.285] and [9, p.272].

Dividing (2.7) by $e^{2\pi x} - 1$ and integrating with respect to $x$ we see that

$$\int_0^\infty \frac{\tan^{-1}(x/u)}{e^{2\pi x} - 1} dx = \int_0^\infty \int_0^\infty \frac{e^{-uy} \sin(xy)}{y(e^{2\pi x} - 1)} dx\, dy = \int_0^\infty \frac{e^{-uy}}{y} dy \int_0^\infty \frac{\sin(xy)}{e^{2\pi x} - 1} dx$$

Since 
$$\frac{1}{e^{2\pi x} - 1} = \sum_{n=1}^\infty e^{-2n\pi x}$$

we have

$$\int_0^\infty \frac{\sin(xy)}{e^{2\pi x} - 1} dx = \sum_{n=1}^\infty \int_0^\infty e^{-2n\pi x} \sin(xy) dx$$

and using (2.5) this becomes

(2.8) $$\int_0^\infty \frac{\sin(xy)}{e^{2\pi x} - 1} dx = \sum_{n=1}^\infty \frac{y}{y^2 + (2n\pi)^2}$$

We then use the identity [16, p.296]

(2.9) $$2\sum_{n=1}^\infty \frac{y}{y^2 + (2n\pi)^2} = \frac{1}{e^y - 1} - \frac{1}{y} + \frac{1}{2}$$

to determine Legendre's relation [49, p.122]

(2.10) $$2\int_0^\infty \frac{\sin(xy)}{e^{2\pi x} - 1} dx = \frac{1}{e^y - 1} - \frac{1}{y} + \frac{1}{2} = \frac{1}{2} \coth\frac{y}{2} - \frac{1}{y}$$

A rigorous derivation of this result is shown in Bromwich's book [16, p.501]. Hence we obtain

$$(2.11) \qquad 2\int_0^\infty \frac{\tan^{-1}(x/u)}{e^{2\pi x}-1}dx = \int_0^\infty \frac{e^{-uy}}{y}\left[\frac{1}{e^y-1}-\frac{1}{y}+\frac{1}{2}\right]dy$$

Then using (2.4) we end up with Binet's first formula for $\log\Gamma(u)$ [49, p.249]

$$(2.12) \qquad \log\Gamma(u) = \left(u-\frac{1}{2}\right)\log u - u + \frac{1}{2}\log(2\pi) + \int_0^\infty \frac{e^{-uy}}{y}\left[\frac{1}{e^y-1}-\frac{1}{y}+\frac{1}{2}\right]dy$$

In passing, we may note that letting $u \to \infty$ in (2.12) enables us to formally derive Stirling's asymptotic approximation for $\log\Gamma(u)$

$$\log\Gamma(u) \approx \left(u-\frac{1}{2}\right)\log u - u + \frac{1}{2}\log(2\pi)$$

Differentiating (2.4) results in

$$(2.13) \qquad \psi(u) = -\frac{1}{2u} + \log u - 2\int_0^\infty \frac{x}{(u^2+x^2)(e^{2\pi x}-1)}dx$$

and since $\psi(u) = -\gamma_0(u)$ we see that this is equivalent to letting $n = 0$ in (1.5). Equation (2.13) was also reported in [49, p.251].

Differentiating (2.12) gives us

$$(2.14) \qquad \psi(u) = -\frac{1}{2u} + \log u + \int_0^\infty e^{-uy}\left[\frac{1}{e^y-1}-\frac{1}{y}+\frac{1}{2}\right]dy$$

As shown in Bromwich's book [16, p.505] we may multiply (2.10) by $e^{-uy}$ and integrate with respect to $y$ to obtain

$$(2.15) \qquad 2\int_0^\infty \frac{x}{(u^2+x^2)(e^{2\pi x}-1)}dx = \int_0^\infty e^{-uy}\left[\frac{1}{e^y-1}-\frac{1}{y}+\frac{1}{2}\right]dy$$

which conforms with (2.13) and (2.14).

Integrating (2.15) over the interval $[0,\infty)$ we again obtain

$$(2.16) \qquad 2\int_0^\infty \frac{\tan^{-1}(x/u)}{e^{2\pi x}-1}dx = \int_0^\infty \frac{e^{-ut}}{t}\left[\frac{1}{e^t-1}-\frac{1}{t}+\frac{1}{2}\right]dt$$

We now integrate (2.16) with respect to $u$. Reversing the order of integration

$$\int_0^v du \int_0^\infty \frac{\tan^{-1}(x/u)}{e^{2\pi x}-1} dx = \int_0^\infty \frac{1}{e^{2\pi x}-1} dx \int_0^t \tan^{-1}(x/u) du$$

and using integration by parts we have

$$\int \tan^{-1}(x/u) du = u\tan^{-1}(x/u) + \frac{1}{2} x \log(u^2 + x^2)$$

Using L'Hôpital's rule we see that $\lim_{u \to 0}\left[u\tan^{-1}(x/u)\right] = 0$ and we therefore obtain

$$\int_0^v \tan^{-1}(x/u) du = v\tan^{-1}(x/v) + \frac{1}{2} x \log(v^2 + x^2) - x \log x$$

We then deduce that

(2.17) $\quad I = 2\int_0^v du \int_0^\infty \frac{\tan^{-1}(x/u)}{e^{2\pi x}-1} dx = 2v\int_0^\infty \frac{\tan^{-1}(x/v)}{e^{2\pi x}-1} dx + \int_0^\infty \frac{x\log(v^2+x^2)}{e^{2\pi x}-1} dx - 2\int_0^\infty \frac{x\log x}{e^{2\pi x}-1} dx$

An evaluation of the following well-known integral is contained for example in [22]

(2.18) $\quad \int_0^\infty \frac{x \log x}{e^{2\pi x}-1} dx = \frac{1}{2}\varsigma'(-1)$

and using Binet's second formula (2.4) for $\log \Gamma(u)$ we obtain

$$I = v\log\Gamma(v) - v\left(v - \frac{1}{2}\right)\log v + v^2 - \frac{1}{2}v\log(2\pi) + \int_0^\infty \frac{x\log(v^2+x^2)}{e^{2\pi x}-1} dx - \varsigma'(-1)$$

Hence we have

(2.19) $\quad \int_0^\infty \frac{x\log(v^2+x^2)}{e^{2\pi x}-1} dx + v\log\Gamma(v) - v\left(v - \frac{1}{2}\right)\log v + v^2 - \frac{1}{2}v\log(2\pi) - \varsigma'(-1)$

$$= \int_0^\infty \frac{1-e^{-vt}}{t^2}\left[\frac{1}{e^t-1} - \frac{1}{t} + \frac{1}{2}\right] dt$$

and with $v = 0$, and noting that $v \log \Gamma(v) = v \log \Gamma(v+1) - v \log v$, we see that (2.18) results.

With $v = 1$ we have

$$\int_0^\infty \frac{x \log(1+x^2)}{e^{2\pi x}-1} dx + 1 - \frac{1}{2}\log(2\pi) - \varsigma'(-1) = \int_0^\infty \frac{1-e^{-t}}{t^2} \left[\frac{1}{e^t-1} - \frac{1}{t} + \frac{1}{2}\right] dt$$

and we shall see in (6.8) that

$$\int_0^\infty \frac{x \log(1+x^2)}{e^{2\pi x}-1} dx = \varsigma'(-1) - \frac{3}{4} + \frac{1}{2}\log(2\pi)$$

Therefore we obtain

(2.20) $$\int_0^\infty \frac{1-e^{-t}}{t^2} \left[\frac{1}{e^t-1} - \frac{1}{t} + \frac{1}{2}\right] dt = \frac{1}{4}$$

Alternatively, we may also integrate (2.19) with respect to $v$ using

$$\int \log(v^2 + x^2) dv = v \log(v^2 + x^2) + 2x \tan^{-1}(v/x) - 2v$$

□

In a similar way to (2.7) we may also show that [37, p.130]

(2.21) $$\frac{1}{2}\log\left(1 + \frac{x^2}{u^2}\right) = \frac{1}{2}\log(u^2 + x^2) - \log u = \int_0^\infty \frac{e^{-uy}[1-\cos(xy)]}{y} dy$$

Multiplying this by $x/(e^{2\pi x} - 1)$ and integrating with respect to $x$ we see that

$$\int_0^\infty \frac{x \log(u^2 + x^2)}{e^{2\pi x}-1} dx - 2\log u \int_0^\infty \frac{x}{e^{2\pi x}-1} dx = 2\int_0^\infty \int_0^\infty \frac{e^{-uy} x[1-\cos(xy)]}{y(e^{2\pi x}-1)} dx\, dy$$

$$= 2\int_0^\infty \frac{e^{-uy}}{y} dy \int_0^\infty \frac{x[1-\cos(xy)]}{e^{2\pi x}-1} dx$$

Employing (6.21) this becomes

$$\int_0^\infty \frac{x\log(u^2+x^2)}{e^{2\pi x}-1}dx - \frac{1}{12}\log u = 2\int_0^\infty \frac{e^{-uy}}{y}dy\int_0^\infty \frac{x[1-\cos(xy)]}{e^{2\pi x}-1}dx$$

Since
$$\frac{1}{e^{2\pi x}-1} = \sum_{n=1}^\infty e^{-2n\pi x}$$

we have

$$\int_0^\infty \frac{x\cos(xy)}{e^{2\pi x}-1}dx = \sum_{n=1}^\infty \int_0^\infty xe^{-2n\pi x}\cos(xy)dx$$

The following integral is fairly elementary

$$\int_0^\infty xe^{-2n\pi x}\cos(xy)dx = \frac{(2n\pi)^2 - y^2}{[(2n\pi)^2 + y^2]^2}$$

and we then have

(2.22) $$\int_0^\infty \frac{\cos(xy)}{e^{2\pi x}-1}dx = \sum_{n=1}^\infty \frac{(2n\pi)^2 - y^2}{[(2n\pi)^2 + y^2]^2}$$

Differentiating (2.9) gives us

$$2\sum_{n=1}^\infty \frac{(2n\pi)^2 - y^2}{[y^2+(2n\pi)^2]^2} = -\frac{e^y}{(e^y-1)^2} + \frac{1}{y^2}$$

and therefore we have

(2.23) $$2\int_0^\infty \frac{\cos(xy)}{e^{2\pi x}-1}dx = -\frac{e^y}{(e^y-1)^2} + \frac{1}{y^2}$$

Hence we obtain

(2.24) $$\int_0^\infty \frac{x\log(u^2+x^2)}{e^{2\pi x}-1}dx - \frac{1}{12}\log u = \int_0^\infty \frac{e^{-uy}}{y}\left[\frac{e^y}{(e^y-1)^2} - \frac{1}{y^2}\right]dy$$

and with $u=1$ we have

(2.25) $$\int_0^\infty \frac{x\log(1+x^2)}{e^{2\pi x}-1}dx = \int_0^\infty \frac{e^{-y}}{y}\left[\frac{e^y}{(e^y-1)^2} - \frac{1}{y^2}\right]dy$$

Using $\log u = \int_0^\infty \frac{e^{-y} - e^{-uy}}{y} dy$ we may write (2.21) as

$$\frac{1}{2}\log(u^2 + x^2) = \int_0^\infty \frac{e^{-uy}[1 - \cos(xy)]}{y} dy + \int_0^\infty \frac{e^{-y} - e^{-uy}}{y} dy$$

$$= \int_0^\infty \frac{e^{-y} - \cos(xy)e^{-uy}}{y} dy$$

$$= \int_0^\infty \frac{e^{-y} - \cos(uy)e^{-xy}}{y} dy \quad \text{(by symmetry)}$$

and we see that

$$\int_0^\infty \frac{e^{-y} - \cos(xy)e^{-uy}}{y} dy = \int_0^\infty \frac{e^{-y} - \cos(uy)e^{-xy}}{y} dy$$

$$\int_0^\infty \frac{\cos(uy)e^{-xy} - \cos(xy)e^{-uy}}{y} dy = 0$$

Differentiation gives us

$$\frac{x}{u^2 + x^2} = \int_0^\infty e^{-uy} \sin(xy) dy = \int_0^\infty e^{-xy} \cos(uy) dy$$

$$\frac{u}{u^2 + x^2} = \int_0^\infty e^{-uy} \cos(xy) dy = \int_0^\infty e^{-xy} \sin(uy) dy$$

□

It may be noted from (2.13) that

(2.26)  $\lim_{u \to \infty}[\psi(u) - \log u] = 0$

and also that

$$\lim_{u \to \infty} u[\psi(u) - \log u] = -\frac{1}{2}$$

We also consider the discrete case with $u = n$

$$\lim_{n \to \infty} n[\psi(n) - \log n] = \lim_{n \to \infty} n\left[\psi(n+1) - \frac{1}{n} - \log n\right]$$

$$= \lim_{n \to \infty} n\left[H_n - \gamma - \frac{1}{n} - \log n\right]$$

$$= 1 + \lim_{n \to \infty} n[H_n - \gamma - \log n]$$

which is in agreement with the more familiar limit

$$\lim_{n \to \infty}\left[n(H_n - \log n - \gamma)\right] = \frac{1}{2}$$

a further proof of which is contained in equation (E.25) of [25].

We also have from (2.13) with $u = n + \frac{1}{2}$

$$\lim_{n \to \infty}\left(n + \frac{1}{2}\right)\left[\psi\left(n + \frac{1}{2}\right) - \log\left(n + \frac{1}{2}\right)\right] = 0$$

We see that

$$\lim_{n \to \infty}\left(n + \frac{1}{2}\right)\left[\psi\left(n + \frac{1}{2}\right) - \log\left(n + \frac{1}{2}\right)\right]$$

$$= \lim_{n \to \infty} n\left[\psi\left(n + \frac{1}{2}\right) - \log\left(n + \frac{1}{2}\right)\right] + \lim_{n \to \infty} \frac{1}{2}\left[\psi\left(n + \frac{1}{2}\right) - \log\left(n + \frac{1}{2}\right)\right]$$

$$= \lim_{n \to \infty} n\left[\psi\left(n + \frac{1}{2}\right) - \log\left(n + \frac{1}{2}\right)\right]$$

and using [46, p.20]

$$\psi\left(n + \frac{1}{2}\right) = -\gamma - 2\log 2 + 2\sum_{k=0}^{n-1} \frac{1}{2k+1}$$

we obtain the known result

$$\lim_{n\to\infty} n\left[2\sum_{k=0}^{n-1}\frac{1}{2k+1}-\gamma-2\log 2-\log\left(n+\frac{1}{2}\right)\right]=0$$

Differentiating (2.13) results in

(2.27) $$\psi'(u)=\frac{1}{2u^2}+\frac{1}{u}+4u\int_0^\infty \frac{x}{(u^2+x^2)^2(e^{2\pi x}-1)}dx$$

which implies that

$$\lim_{u\to\infty}\left[u^2\psi'(u)-u\right]=\frac{1}{2}$$

**3. An application of the Abel-Plana summation formula**

Adamchik [5] noted that the Hermite integral for the Hurwitz zeta function may be derived from the Abel-Plana summation formula [46, p.90]

(3.1) $$\sum_{k=0}^\infty f(k)=\frac{1}{2}f(0)+\int_0^\infty f(x)dx+i\int_0^\infty \frac{f(ix)-f(-ix)}{e^{2\pi x}-1}dx$$

which applies to functions which are analytic in the right-hand plane and satisfy the convergence condition

$$\lim_{y\to\infty} e^{-2\pi y}|f(x+iy)|=0$$

uniformly on any finite interval of $x$. Derivations of the Abel-Plana summation formula are given in [49, p.145] and [50, p.108].

Letting $f(k)=(k+u)^{-s}$ we obtain

(3.2) $$\varsigma(s,u)=\sum_{k=0}^\infty \frac{1}{(k+u)^s}=\frac{u^{-s}}{2}+\frac{u^{1-s}}{s-1}+i\int_0^\infty \frac{(u+ix)^{-s}-(u-ix)^{-s}}{e^{2\pi x}-1}dx$$

Then, noting that

$$(u+ix)^{-s}-(u-ix)^{-s}=(re^{i\theta})^{-s}-(re^{-i\theta})^{-s}$$

$$=r^{-s}[e^{-is\theta}-e^{is\theta}]$$

$$= \frac{2}{i(u^2+x^2)^{s/2}} \sin(s\tan^{-1}(x/u))$$

we may write (3.2) as Hermite's integral for $\varsigma(s,u)$

(3.3) $$\varsigma(s,u) = \frac{u^{-s}}{2} + \frac{u^{1-s}}{s-1} + 2\int_0^\infty \frac{\sin(s\tan^{-1}(x/u))}{(u^2+x^2)^{s/2}(e^{2\pi x}-1)} dx$$

We now take one step back and differentiate the intermediate equation (3.2) with respect to $s$ to obtain

(3.4)
$$\varsigma'(s,u) = -\frac{1}{2}u^{-s}\log u - \frac{u^{1-s}[1+(s-1)\log u]}{(s-1)^2} + i\int_0^\infty \frac{(u-ix)^{-s}\log(u-ix)-(u+ix)^{-s}\log(u+ix)}{e^{2\pi x}-1} dx$$

It will be seen several times in this paper that computationally it is much easier to deal with (3.2) rather than Hermite's integral (3.3).

It may be noted that

$$i\int_0^\infty \frac{(u-ix)^{-s}\log(u-ix)-(u+ix)^{-s}\log(u+ix)}{e^{2\pi x}-1} dx = 2\operatorname{Im}\int_0^\infty \frac{(u+ix)^{-s}\log(u+ix)}{e^{2\pi x}-1} dx$$

With $s=0$ in (3.4) we obtain

$$\varsigma'(0,u) = \left(u-\frac{1}{2}\right)\log u - u + i\int_0^\infty \frac{\log(u-ix)-\log(u+ix)}{e^{2\pi x}-1} dx$$

This may be written as

$$\varsigma'(0,u) = \left(u-\frac{1}{2}\right)\log u - u + 2\int_0^\infty \frac{\tan^{-1}(x/u)}{e^{2\pi x}-1} dx$$

and, using Lerch's identity (5.4), we see that this is equivalent to (2.4).

When $s=-1$ we have

$$\varsigma'(-1,u) = \frac{1}{2}u(u-1)\log u - \frac{1}{4}u^2 + i\int_0^\infty \frac{(u-ix)\log(u-ix)-(u+ix)\log(u+ix)}{e^{2\pi x}-1} dx$$

or

$$\varsigma'(-1,u) = \frac{1}{2}u(u-1)\log u - \frac{1}{4}u^2 + \int_0^\infty \frac{x\log(u^2+x^2) + 2u\tan^{-1}(x/u)}{e^{2\pi x}-1}dx$$

and we will see that this is equivalent to (6.11). Using (2.4) we obtain

(3.5) $$\int_0^\infty \frac{x\log(u^2+x^2)}{e^{2\pi x}-1}dx = \varsigma'(-1,u) - \frac{1}{2}u^2\log u - \frac{3}{4}u^2 + \frac{1}{2}u\log(2\pi) - u\log\Gamma(u)$$

which was previously derived by Adamchik [3].

In particular, with $u=0$ we have

$$\varsigma'(-1) = \varsigma'(-1,0) = 2\int_0^\infty \frac{x\log x}{e^{2\pi x}-1}dx$$

which we noted previously in (2.18).

With $u=1$ we get

$$\varsigma'(-1) = \varsigma'(-1,1) = -\frac{1}{4} + \int_0^\infty \frac{x\log(1+x^2)}{e^{2\pi x}-1}dx + 2\int_0^\infty \frac{\tan^{-1}x}{e^{2\pi x}-1}dx$$

and from (2.14) we see that

$$\int_0^\infty \frac{\tan^{-1}x}{e^{2\pi x}-1}dx = \frac{1}{2} - \frac{1}{4}\log(2\pi)$$

We then determine that

(3.6) $$\int_0^\infty \frac{x\log(1+x^2)}{e^{2\pi x}-1}dx = \varsigma'(-1) - \frac{3}{4} + \frac{1}{2}\log(2\pi)$$

Similarly, letting $s=-2$ gives us

$$\varsigma'(-2,u) = -\frac{1}{2}u^2\log u - \frac{1}{9}u^3(1-3\log u) + i\int_0^\infty \frac{(u-ix)^2\log(u-ix) - (u+ix)^2\log(u+ix)}{e^{2\pi x}-1}dx$$

$$= -\frac{1}{2}u^2\log u - \frac{1}{9}u^3(1-3\log u) + 2\int_0^\infty \frac{(u^2-x^2)\tan^{-1}(x/u) + ux\log(u^2+x^2)}{e^{2\pi x}-1}dx$$

$$= -\frac{1}{2}u^2 \log u - \frac{1}{9}u^3(1-3\log u) + 2u^2 \int_0^\infty \frac{\tan^{-1}(x/u)}{e^{2\pi x}-1}dx - 2\int_0^\infty \frac{x^2 \tan^{-1}(x/u)}{e^{2\pi x}-1}dx$$

$$+2u\int_0^\infty \frac{x\log(u^2+x^2)}{e^{2\pi x}-1}dx$$

and we may then use (2.4) and (3.5) for the first and third integrals respectively. This gives us

(3.7)

$$2\int_0^\infty \frac{x^2 \tan^{-1}(x/u)}{e^{2\pi x}-1}dx = -\frac{1}{2}u^2 \log u - \frac{1}{9}u^3(1-3\log u) + u^2\left(\log\Gamma(u) - \left(u-\frac{1}{2}\right)\log u + u - \frac{1}{2}\log(2\pi)\right)$$

$$+2u\left(\varsigma'(-1,u) - \frac{1}{2}u^2 \log u - \frac{3}{4}u^2 + \frac{1}{2}u\log(2\pi) - u\log\Gamma(u)\right) - \varsigma'(-2,u)$$

which was also previously evaluated by Adamchik [3].

We also note that

$$\varsigma'(2,u) = -\frac{1}{2u^2}\log u - \frac{1+\log u}{u} + i\int_0^\infty \frac{(u-ix)^{-2}\log(u-ix) - (u+ix)^{-2}\log(u+ix)}{e^{2\pi x}-1}dx$$

$$= -\frac{1}{2u^2}\log u - \frac{1+\log u}{u} + 2\,\text{Im}\int_0^\infty \frac{\log(u+ix)}{(u+ix)^2(e^{2\pi x}-1)}dx$$

$$= -\frac{1}{2u^2}\log u - \frac{1+\log u}{u} + 2\,\text{Im}\int_0^\infty \frac{(u-ix)^2 \log(u+ix)}{(u^2+x^2)^2(e^{2\pi x}-1)}dx$$

Hence we have

$$\varsigma'(2,u) = -\frac{1}{u} - \frac{\log u}{2u^2} - \frac{\log u}{u} + 2\int_0^\infty \frac{(u^2-x^2)\tan^{-1}(x/u) - ux\log(u^2+x^2)}{(u^2+x^2)^2(e^{2\pi x}-1)}dx$$

Differentiating (5.1) gives us

$$\gamma_1'(u) = \frac{1-\log u}{2u^2} - \frac{\log u}{u} + 2\int_0^\infty \frac{ux - ux\log(u^2+x^2)}{(u^2+x^2)^2(e^{2\pi x}-1)}dx + 2u\int_0^\infty \frac{2u\tan^{-1}(x/u)+x}{(u^2+x^2)^2(e^{2\pi x}-1)}dx$$

$$-2\int_0^\infty \frac{\tan^{-1}(x/u)}{(u^2+x^2)(e^{2\pi x}-1)}dx$$

$$= \frac{1-\log u}{2u^2} - \frac{\log u}{u} - 2u\int_0^\infty \frac{x\log(u^2+x^2)}{(u^2+x^2)^2(e^{2\pi x}-1)}dx + 4u^2\int_0^\infty \frac{\tan^{-1}(x/u)}{(u^2+x^2)^2(e^{2\pi x}-1)}dx$$

$$-2\int_0^\infty \frac{\tan^{-1}(x/u)}{(u^2+x^2)(e^{2\pi x}-1)}dx + 4u\int_0^\infty \frac{x}{(u^2+x^2)^2(e^{2\pi x}-1)}dx$$

We recall (2.13)

$$\psi(u) = -\frac{1}{2u} + \log u - 2\int_0^\infty \frac{x}{(u^2+x^2)(e^{2\pi x}-1)}dx$$

and differentiation results in

$$\psi'(u) = \frac{1}{2u^2} + \frac{1}{u} + 4u\int_0^\infty \frac{x}{(u^2+x^2)^2(e^{2\pi x}-1)}dx$$

Substituting the above integral, we therefore obtain

$$\gamma_1'(u) = -\frac{1}{u} - \frac{\log u}{2u^2} - \frac{\log u}{u} - 2u\int_0^\infty \frac{x\log(u^2+x^2)}{(u^2+x^2)^2(e^{2\pi x}-1)}dx + 4u^2\int_0^\infty \frac{\tan^{-1}(x/u)}{(u^2+x^2)^2(e^{2\pi x}-1)}dx$$

$$-2\int_0^\infty \frac{\tan^{-1}(x/u)}{(u^2+x^2)(e^{2\pi x}-1)}dx + \psi'(u)$$

Comparing this with the above result for $\varsigma'(2,u)$ we determine, after a little algebra, that

(3.8) $\qquad \gamma_1'(u) = \varsigma(2,u) + \varsigma'(2,u) = \psi'(u) + \varsigma'(2,u)$

which was previously noted in equation (4.3.223b) in [22].

With $s=2$ in (3.2) we get

(3.9) $\qquad \varsigma(2,u) = \sum_{k=0}^\infty \frac{1}{(k+u)^2} = \frac{1}{2u^2} + \frac{1}{u} + 4u\int_0^\infty \frac{x}{(u^2+x^2)^2(e^{2\pi x}-1)}dx$

and we immediately note that the right-hand side is the same as (2.21), leading to the well-known result $\varsigma(2,u) = \psi'(u)$.

We also see from (3.2) that when $s$ is a negative integer

$$\varsigma(-m,u) = \frac{u^m}{2} - \frac{u^{m+1}}{m+1} + i\int_0^\infty \frac{(u+ix)^m - (u-ix)^m}{e^{2\pi x} - 1} dx$$

$$= \frac{u^m}{2} - \frac{u^{m+1}}{m+1} - 2\,\text{Im}\int_0^\infty \frac{(u+ix)^m}{e^{2\pi x} - 1} dx$$

Using the following expression for the Bernoulli polynomials

$$B_m(u) = \sum_{k=0}^m \binom{m}{k} B_k u^{m-k}$$

we easily find the well-known formula [7, p.264]

(3.10) $$\varsigma(-m,u) = -\frac{B_{m+1}(u)}{m+1}$$

This formula was also derived by Boros et al [12], albeit in a more convoluted manner using Hermite's integral (this is another example of where (3.2) is easier to manipulate than (3.3)).

We now consider the second derivative of the Hurwitz zeta function

$$\varsigma''(s,u) = \frac{1}{2} u^{-s} \log^2 u + \frac{u^{1-s}[(s-1)\log u + 2]}{(s-1)^3} + \frac{u^{1-s} \log^2 u}{s-1} + \frac{u^{1-s} \log u}{(s-1)^2}$$

$$-i\int_0^\infty \frac{(u-ix)^{-s} \log^2(u-ix) - (u+ix)^{-s} \log^2(u+ix)}{e^{2\pi x} - 1} dx$$

where with $s = 0$ we have

(3.11) $$\varsigma''(0,u) = \left(\frac{1}{2} - u\right) \log^2 u + 2u \log u - 2u - 2\int_0^\infty \frac{\log(u^2 + x^2) \tan^{-1}(x/u)}{e^{2\pi x} - 1} dx$$

and we shall see this again in (5.3).

We note in passing that with the change of variables $x = uv$ we have

$$\int_0^\infty \frac{x \log(u^2+x^2)}{(u^2+x^2)(e^{2\pi x}-1)} dx = \int_0^\infty \frac{u^2 v \log(u^2+(uv)^2)}{(u^2+(uv)^2)(e^{2\pi uv}-1)} dv$$

$$= \int_0^\infty \frac{v \log(1+v^2)+2uv \log u}{(1+v^2)(e^{2\pi uv}-1)} dv$$

$$= \int_0^\infty \frac{v \log(1+v^2)}{(1+v^2)(e^{2\pi uv}-1)} dv + 2u \log u \int_0^\infty \frac{v}{(1+v^2)(e^{2\pi uv}-1)} dv$$

**4. Another derivation of Coffey's integral for the Stieltjes constants**

A slight modification of (3.2) results in

(4.1) $$\varsigma(s,u) - \frac{1}{s-1} = \frac{u^{-s}}{2} + \frac{u^{1-s}}{s-1} - \frac{1}{s-1} + i \int_0^\infty \frac{(u+ix)^{-s}-(u-ix)^{-s}}{e^{2\pi x}-1} dx$$

It may be seen from (1.6) that

$$\gamma_n(u) = (-1)^n \lim_{s \to 1} \frac{d^n}{ds^n}\left[\varsigma(s,u) - \frac{1}{s-1}\right] = (-1)^n \lim_{s \to 1}\left[\varsigma^{(n)}(s,u) - \frac{(-1)^n n!}{(s-1)^{n+1}}\right]$$

and differentiating (4.1) $n$ times gives us

$$\frac{d^n}{ds^n}\left[\varsigma(s,u) - \frac{1}{s-1}\right]$$

$$= \frac{u^{-s}(-1)^n \log^n u}{2} + f^{(n)}(s) + i(-1)^n \int_0^\infty \frac{(u+ix)^{-s}\log^n(u+ix)-(u-ix)^{-s}\log^n(u-ix)}{e^{2\pi x}-1} dx$$

where we have denoted $f(s)$ as

$$f(s) = \frac{u^{1-s}-1}{s-1}$$

We can represent $f(s)$ by the following integral

$$f(s) = \frac{u^{1-s}-1}{s-1} = -\int_1^u x^{-s} dx$$

so that

$$f^{(n)}(s) = -(-1)^n \int_1^u x^{-s} \log^n dx$$

and thus

(4.2) $$f^{(n)}(1) = -(-1)^n \int_1^u \frac{\log^n x}{x} dx = (-1)^{n+1} \frac{\log^{n+1} u}{n+1}$$

Therefore, upon taking the limit $s \to 1$, we obtain

(4.3) $$\gamma_n(u) = \frac{1}{2u} \log^n u - \frac{1}{n+1} \log^{n+1} u + i \int_0^\infty \frac{(u-ix)\log^n(u+ix) - (u+ix)\log^n(u-ix)}{(u^2+x^2)(e^{2\pi x}-1)} dx$$

which is equivalent to Coffey's formula (1.5). This derivation is slightly more direct than the one originally provided by Coffey [20] in 2007.

$\square$

An alternative proof of (4.2) is shown below. Let us now define $f_1(s)$ more generally by

$$f_1(s) = \frac{g(s,u)-1}{s-1}$$

where $g(1,u) = 1$ and $g^{(n)}(s,u)$ exists.

We easily see that

$$\frac{d^{n+1}}{ds^{n+1}}[(s-1)f_1(s)] = g^{(n+1)}(s,u)$$

The following result is readily derived using the Leibniz differentiation formula

$$\frac{d^{n+1}}{ds^{n+1}}[(s-1)f_1(s)] = (s-1)f_1^{(n+1)}(s) + (n+1) f_1^{(n)}(s)$$

and hence we see that

(4.4) $$(s-1)f_1^{(n+1)}(s) + (n+1) f_1^{(n)}(s) = g^{(n+1)}(s,u)$$

With $s=1$ in (4.4) we have (as shown below)

(4.5) $$(n+1)f_1^{(n)}(1) = g^{(n+1)}(1,u)$$

We now prove that (4.5) is valid. Differentiation results in

$$f_1'(s) = \frac{(s-1)g'(s,u) - g(s,u) + 1}{(s-1)^2}$$

and we note that we have the indeterminate form $f_1'(1) \approx \frac{0}{0}$. Therefore, applying L'Hôpital's rule we see that

$$\lim_{s \to 1} f'(s) = \lim_{s \to 1} \frac{(s-1)g''(s,u)}{2(s-1)} = \frac{1}{2}g''(1,u)$$

Hence we see that $f_1^{(1)}(1)$ exists and from (4.4) we deduce that $f_1^{(2)}(1)$ also exists. Therefore, by induction we see that $f_1^{(n+1)}(1)$ also exists and we then determine that $\lim_{s \to 1}(s-1)f_1^{(n+1)}(s) = 0$, thereby proving (4.5) above. It is clear that in the case where $g(s,u) = u^{1-s}$ then $g^{(n+1)}(s,u) = u^{1-s}(-1)^{n+1}\log^{n+1} u$. [The logic employed above is rather incomplete and requires further clarification].

The above analysis may indeed be employed to make a more succinct derivation of the formula (2.1) for the Stieltjes constants. We start with the Hasse identity [31] for the Hurwitz zeta function which is valid for all $s \in \mathbf{C}$ provided $\mathrm{Re}(s) \neq 1$

(4.6) $$(s-1)\varsigma(s,u) = \sum_{i=0}^{\infty} \frac{1}{i+1} \sum_{j=0}^{i} \binom{i}{j} \frac{(-1)^j}{(u+j)^{s-1}}$$

We have the limit

(4.7) $$\lim_{s \to 1}(s-1)\varsigma(s,u) = \sum_{i=0}^{\infty} \frac{1}{i+1} \sum_{j=0}^{i} \binom{i}{j}(-1)^j = \sum_{i=0}^{\infty} \frac{1}{i+1}\delta_{i,0} = 1$$

We now let $g(s,u) = (s-1)\varsigma(s,u)$ and consider

$$f_2(s) = \frac{(s-1)\varsigma(s,u) - 1}{s-1} = \varsigma(s,u) - \frac{1}{s-1}$$

and from (4.7) we see that $g(1,u) = 1$.

Hence, using (4.4) we have

$$f_2^{(n)}(1) = \frac{1}{n+1}g^{(n+1)}(1,u) = \frac{(-1)^{n+1}}{n+1}\sum_{i=0}^{\infty} \frac{1}{i+1}\sum_{j=0}^{i} \binom{i}{j}(-1)^j \log^{n+1}(u+j)$$

Hence we immediately get another proof of (2.1)

$$(4.8) \qquad \gamma_n(u) = -\frac{1}{n+1}\sum_{i=0}^{\infty}\frac{1}{i+1}\sum_{j=0}^{i}\binom{i}{j}(-1)^j \log^{n+1}(u+j)$$

For future reference we note from (4.6) that

$$(4.9) \qquad \frac{d^{n+1}}{ds^{n+1}}[(s-1)\varsigma(s,u)] = (-1)^{n+1}\sum_{i=0}^{\infty}\frac{1}{i+1}\sum_{j=0}^{i}\binom{i}{j}\frac{(-1)^j \log^{n+1}(u+j)}{(u+j)^{s-1}}$$

and thus

$$\left.\frac{d^{n+1}}{ds^{n+1}}[(s-1)\varsigma(s,u)]\right|_{s=1} = (-1)^{n+1}\sum_{i=0}^{\infty}\frac{1}{i+1}\sum_{j=0}^{i}\binom{i}{j}(-1)^j \log^{n+1}(u+j)$$

Reference to (4.8) allows us to conclude that for $n \geq 0$

$$(4.10) \qquad \left.\frac{d^{n+1}}{ds^{n+1}}[(s-1)\varsigma(s,u)]\right|_{s=1} = (-1)^n (n+1)\gamma_n(u)$$

□

We now refer to the Hasse/Sondow identity (see [31] and [45])

$$(4.11) \qquad \varsigma_a(s) = \sum_{i=0}^{\infty}\frac{1}{2^{i+1}}\sum_{j=0}^{i}\binom{i}{j}\frac{(-1)^j}{(1+j)^s}$$

where the alternating Riemann zeta function $\varsigma_a(s)$ is defined by

$$\varsigma_a(s) = (1-2^{1-s})\varsigma(s)$$

We then see that

$$(4.12) \qquad (s-1)\varsigma(s) = \frac{s-1}{1-2^{1-s}}\sum_{i=0}^{\infty}\frac{1}{2^{i+1}}\sum_{j=0}^{i}\binom{i}{j}\frac{(-1)^j}{(1+j)^s}$$

We define $f_3(s)$ by

$$f_3(s) = \frac{2^{1-s}-1}{s-1}$$

and hence, having regard to (4.5), we have

$$(4.13) \qquad f_3^{(n)}(1) = \frac{\log^{n+1} 2}{n+1}$$

We see that

$$\varsigma_a(s) = \frac{1-2^{1-s}}{s-1}(s-1)\varsigma(s)$$

and hence

$$\lim_{s\to 1}\varsigma_a(s) = \lim_{s\to 1}\frac{1-2^{1-s}}{s-1}\lim_{s\to 1}[(s-1)\varsigma(s)]$$

Using L'Hôpital's rule we have

$$\lim_{s\to 1}\frac{1-2^{1-s}}{s-1} = \log 2 = -f_3(1)$$

and, combining this with (4.6), we note the well-known limit

(4.14) $$\lim_{s\to 1}\varsigma_a(s) = \log 2$$

We now differentiate (4.12) to obtain

(4.15)

$$\frac{d}{ds}[(s-1)\varsigma(s)] = -\frac{d}{ds}\left[\frac{1}{f_3(s)}\right]\sum_{i=0}^{\infty}\frac{1}{2^{i+1}}\sum_{j=0}^{i}\binom{i}{j}\frac{(-1)^j}{(1+j)^s} - \frac{s-1}{1-2^{1-s}}\sum_{i=0}^{\infty}\frac{1}{2^{i+1}}\sum_{j=0}^{i}\binom{i}{j}\frac{(-1)^j\log(1+j)}{(1+j)^s}$$

and with $s=1$ we have using (4.10) and (4.13)

(4.16) $$\gamma_0 = \gamma = \frac{1}{2}\log 2 - \frac{1}{\log 2}\sum_{i=0}^{\infty}\frac{1}{2^{i+1}}\sum_{j=0}^{i}\binom{i}{j}(-1)^j\frac{\log(1+j)}{1+j}$$

This expression for Euler's constant was originally derived by Coffey [19] in 2006 (a different derivation is contained in equation (4.4.116g) in [23])).

Equation (4.16) may be written as

(4.16.1) $$\gamma = \frac{1}{2}\log 2 - \frac{1}{\log 2}\varsigma_a^{(1)}(1)$$

Differentiating (4.15) results in

$$\frac{d^2}{ds^2}[(s-1)\varsigma(s)] = \sum_{i=0}^{\infty}\frac{1}{2^{i+1}}\sum_{j=0}^{i}\binom{i}{j}\frac{(-1)^j \log(1+j)}{(1+j)^s}\frac{d}{ds}\left[\frac{1}{f_3(s)}\right] - \sum_{i=0}^{\infty}\frac{1}{2^{i+1}}\sum_{j=0}^{i}\binom{i}{j}\frac{(-1)^j}{(1+j)^s}\frac{d^2}{ds^2}\left[\frac{1}{f_3(s)}\right]$$

$$-\frac{1}{f_3(s)}\sum_{i=0}^{\infty}\frac{1}{2^{i+1}}\sum_{j=0}^{i}\binom{i}{j}\frac{(-1)^j \log^2(1+j)}{(1+j)^s} + \sum_{i=0}^{\infty}\frac{1}{2^{i+1}}\sum_{j=0}^{i}\binom{i}{j}\frac{(-1)^j \log(1+j)}{(1+j)^s}\frac{d}{ds}\left[\frac{1}{f_3(s)}\right]$$

and with $s=1$ we deduce that the second Stieltjes constant may be represented by

(4.16.2)
$$\gamma_1 = -\frac{1}{12}\log^2 2 + \frac{1}{2}\sum_{i=0}^{\infty}\frac{1}{2^{i+1}}\sum_{j=0}^{i}\binom{i}{j}(-1)^j\frac{\log(1+j)}{1+j} - \frac{1}{2\log 2}\sum_{i=0}^{\infty}\frac{1}{2^{i+1}}\sum_{j=0}^{i}\binom{i}{j}(-1)^j\frac{\log^2(1+j)}{1+j}$$

which was also previously determined by Coffey [19] in a slightly different manner.

Equation (4.16.2) may be written as

(4.16.3) $$\gamma_1 = -\frac{1}{12}\log^2 2 - \frac{1}{2}\varsigma_a^{(1)}(1) - \frac{1}{2\log 2}\varsigma_a^{(2)}(1)$$

We also have the generalisation (see equation (4.4.79) in [23])

(4.17) $$\varsigma_a(s,u) = \sum_{i=0}^{\infty}\frac{1}{2^{i+1}}\sum_{j=0}^{i}\binom{i}{j}\frac{(-1)^j}{(u+j)^s}$$

where $\varsigma_a(s,u)$ may be regarded as an alternating Hurwitz zeta function which may be written as

(4.18) $$\varsigma_a(s,u) = \sum_{i=0}^{\infty}\frac{(-1)^i}{(i+u)^s}$$

This would enable us to determine similar expressions for $\gamma_n(u)$.

We note from G&R [29, p.897] that

(4.19) $$\varsigma_a(1,u) = \sum_{i=0}^{\infty}\frac{(-1)^i}{i+u} = \frac{1}{2}\left[\psi\left(\frac{1+u}{2}\right) - \psi\left(\frac{u}{2}\right)\right]$$

and this gives us

$$\varsigma_a(1,1) = \varsigma_a(1) = \sum_{i=0}^{\infty}\frac{(-1)^i}{i+1} = \frac{1}{2}\left[\psi(1) - \psi\left(\frac{1}{2}\right)\right]$$

and since [46, p.20]

$$\psi\left(\frac{1}{2}\right) = -\gamma - 2\log 2$$

we again see that

$$\varsigma_a(1,1) = \varsigma_a(1) = \log 2$$

We may also consider the expression

$$(1-2^{1-s})(s-1)\varsigma(s) = (s-1)\sum_{i=0}^{\infty}\frac{1}{2^{i+1}}\sum_{j=0}^{i}\binom{i}{j}\frac{(-1)^j}{(1+j)^s}$$

and differentiate this $n+1$ times to obtain for the left-hand side

$$\frac{d^{n+1}}{ds^{n+1}}[(1-2^{1-s})(s-1)\varsigma(s)] = \sum_{k=0}^{n+1}\binom{n+1}{k}\frac{d^k}{ds^k}[(s-1)\varsigma(s)]\frac{d^{n+1-k}}{ds^{n+1-k}}(1-2^{1-s})$$

Using (4.10) to evaluate this at $s=1$ we have

$$\sum_{k=0}^{n+1}\binom{n+1}{k}\frac{d^k}{ds^k}[(s-1)\varsigma(s)]\frac{d^{n+1-k}}{ds^{n+1-k}}(1-2^{1-s})\bigg|_{s=1}$$

$$= (-1)^n \log^{n+1} 2 + \sum_{k=1}^{n}\binom{n+1}{k}(-1)^{k-1}k\gamma_{k-1}(-1)^{n-k}\log^{n+1-k} 2$$

since the $(n+1)$th term vanishes at $s=1$ and we have isolated the first term.

For the right-hand side we have using (4.3)

$$\frac{d^{n+1}}{ds^{n+1}}\left[(s-1)\sum_{i=0}^{\infty}\frac{1}{2^{i+1}}\sum_{j=0}^{i}\binom{i}{j}\frac{(-1)^j}{(1+j)^s}\right]\bigg|_{s=1} = (-1)^n(n+1)\sum_{i=0}^{\infty}\frac{1}{2^{i+1}}\sum_{j=0}^{i}\binom{i}{j}\frac{(-1)^j \log^n(1+j)}{1+j}$$

and therefore we obtain

$$(n+1)\sum_{i=0}^{\infty}\frac{1}{2^{i+1}}\sum_{j=0}^{i}\binom{i}{j}\frac{(-1)^j \log^n(1+j)}{1+j} = \log^{n+1} 2 - \sum_{k=1}^{n}\binom{n+1}{k}k\gamma_{k-1}\log^{n+1-k} 2$$

This may also be written as

$$(-1)^{n+1}\varsigma_a^{(n)}(1) = \frac{1}{n+1}\sum_{k=1}^{n}\binom{n+1}{k}k\gamma_{k-1}\log^{n+1-k}2 - \frac{1}{n+1}\log^{n+1}2$$

and since $\frac{k}{n+1}\binom{n+1}{k} = \binom{n}{k-1}$ we see that

(4.20) $$(-1)^{n+1}\varsigma_a^{(n)}(1) = \sum_{m=1}^{\infty}(-1)^m \frac{\log^n m}{m} = \sum_{k=0}^{n-1}\binom{n}{k}\gamma_k \log^{n-k}2 - \frac{1}{n+1}\log^{n+1}2$$

as shown by Dilcher in [26].

The above formula (4.20) was in fact first reported by Briggs and Chowla [15] in 1955 where they showed that

(4.21) $$\varsigma_a^{(k)}(1) = k!\sum_{r=1}^{k+1}\frac{(-1)^{r+1}\log^r 2}{r!}A_{k-r}$$

with $A_n = \frac{(-1)^n}{n!}\gamma_n$ and $A_{-1} = 1$.

It was rediscovered in 1972 by Liang and Todd [36] who used it to numerically compute the values of the first 20 Stieltjes constants (see Dilcher's paper [26] for more details).

Equation (4.20) was posed as a problem by Klamkin [34] in 1954 for the case $n=1$ and is closely related to an earlier problem posed by Sandham [43] in 1950.

Using simple algebra, equation (4.20) may be written as

$$\frac{1}{\log^n 2}\left[(-1)^{n+1}\varsigma_a^{(n)}(1) + \frac{1}{n+1}\log^{n+1}2 + \gamma_n\right] = \sum_{k=0}^{n}\binom{n}{k}(-1)^k \frac{(-1)^k \gamma_k}{\log^k 2}$$

and applying the binomial inversion formula

$$a_n = \sum_{k=0}^{n}\binom{n}{k}(-1)^k b_k \quad \Leftrightarrow \quad b_n = \sum_{k=0}^{n}\binom{n}{k}(-1)^k a_k$$

we obtain

(4.22) $$\gamma_n = (-1)^n\sum_{k=0}^{n}\binom{n}{k}(-1)^k \log^{n-k}2\left[(-1)^{k+1}\varsigma_a^{(k)}(1) + \frac{1}{k+1}\log^{k+1}2 + \gamma_k\right]$$

This formula immediately gives us (4.16.1) and (4.16.3).

Since $\sum_{k=0}^{n}\binom{n}{k}\frac{(-1)^k}{k+1}=\frac{1}{n+1}$ this may also be written as

(4.23) $\gamma_n = (-1)^n \sum_{k=0}^{n}\binom{n}{k}(-1)^k \log^{n-k} 2\left[(-1)^{k+1}\varsigma_a^{(k)}(1)+\gamma_k\right]+\frac{1}{n+1}(-1)^n \log^{n+1} 2$

## 5. Some aspects of the Stieltjes constants

Letting $n=1$ in (1.5) gives us

$$\gamma_1(u)=\frac{1}{2u}\log u - \frac{1}{2}\log^2 u - 2\operatorname{Re}\int_0^\infty \frac{i(u+ix)\log(u-ix)}{(u^2+x^2)(e^{2\pi x}-1)}dx$$

We easily see that

$$\operatorname{Re}\int_0^\infty \frac{i(u+ix)\log(u-ix)}{(u^2+x^2)(e^{2\pi x}-1)}dx = \operatorname{Re}\int_0^\infty \frac{i(u+ix)\left[\log\sqrt{u^2+x^2}-i\tan^{-1}(x/u)\right]}{(u^2+x^2)(e^{2\pi x}-1)}dx$$

$$=-\frac{1}{2}u\int_0^\infty \frac{x\log(u^2+x^2)}{(u^2+x^2)(e^{2\pi x}-1)}dx + u^2\int_0^\infty \frac{\tan^{-1}(x/u)}{(u^2+x^2)(e^{2\pi x}-1)}dx$$

and we therefore obtain a particular case for the Stieltjes constant $\gamma_1(u)$

(5.1) $\gamma_1(u)=\frac{1}{2u}\log u - \frac{1}{2}\log^2 u + \int_0^\infty \frac{x\log(u^2+x^2)}{(u^2+x^2)(e^{2\pi x}-1)}dx - 2u\int_0^\infty \frac{\tan^{-1}(x/u)}{(u^2+x^2)(e^{2\pi x}-1)}dx$

and using (2.1) this may be expressed as

$$-\frac{1}{2}\sum_{i=0}^{\infty}\frac{1}{i+1}\sum_{j=0}^{i}\binom{i}{j}(-1)^j \log^2(u+j) = \frac{1}{2u}\log u - \frac{1}{2}\log^2 u + \int_0^\infty \frac{x\log(u^2+x^2)}{(u^2+x^2)(e^{2\pi x}-1)}dx$$

$$-2u\int_0^\infty \frac{\tan^{-1}(x/u)}{(u^2+x^2)(e^{2\pi x}-1)}dx$$

It may be noted that Shail [44, p.799] made reference to these "seemingly intractable" integrals in 2000.

Letting $x=ut$ in the integrals results in

$$\gamma_1(u) = \frac{1}{2u}\log u - \frac{1}{2}\log^2 u + 2\log u \int_0^\infty \frac{t}{(1+t^2)(e^{2\pi ut}-1)}\,dt + \int_0^\infty \frac{t\log(1+t^2)}{(1+t^2)(e^{2\pi ut}-1)}\,dt$$

$$-2\int_0^\infty \frac{\tan^{-1} t}{(1+t^2)(e^{2\pi ut}-1)}\,dt$$

Boros and Moll [13] noted that for any differentiable function we have

$$u\frac{\partial}{\partial u} f(ut) = t\frac{\partial}{\partial t} f(ut)$$

and it seem to me that this relationship could be usefully applied to some of the integrals reported in this paper.

It is readily seen that

$$-\frac{d}{du}\left[\log(u^2+x^2)\tan^{-1}(x/u)\right] = \frac{x\log(u^2+x^2) - 2u\tan^{-1}(x/u)}{(u^2+x^2)}$$

and integrating (5.1) therefore results in

$$\int_1^t \gamma_1(u)\,du = \frac{1}{4}\log^2 t + 1 - \frac{1}{2}t(\log^2 t - 2\log t + 2) - \int_0^\infty \frac{\log(t^2+x^2)\tan^{-1}(x/u) - \log(1+x^2)\tan^{-1}(x)}{e^{2\pi x}-1}\,dx$$

We note from [22, Eq. (4.3.231)] that for $t > 0$

$$\int_1^t \gamma_n(u)\,du = \frac{(-1)^{n+1}}{n+1}\left[\varsigma^{(n+1)}(0,t) - \varsigma^{(n+1)}(0)\right]$$

and for $n=1$ we have

$$\int_1^t \gamma_1(u)\,du = \frac{1}{2}\left[\varsigma''(0,t) - \varsigma''(0)\right]$$

We therefore deduce that

(5.2) $$\int_0^\infty \frac{\log(t^2+x^2)\tan^{-1}(x/t) - \log(1+x^2)\tan^{-1}(x)}{e^{2\pi x}-1}\,dx$$

$$= \frac{1}{4}\log^2 t + 1 - \frac{1}{2}t(\log^2 t - 2\log t + 2) - \frac{1}{2}\left[\varsigma''(0,t) - \varsigma''(0)\right]$$

We may also write this as

(5.3) $$\int_0^\infty \frac{\log(t^2+x^2)\tan^{-1}(x/t)}{e^{2\pi x}-1}dx = \frac{1}{4}\log^2 t - \frac{1}{2}t(\log^2 t - 2\log t + 2) - \frac{1}{2}\varsigma''(0,t)$$

which we have seen previously in (3.11).

Differentiation of (5.3) results in

$$-\int_0^\infty \frac{x\log(t^2+x^2)}{(t^2+x^2)(e^{2\pi x}-1)}dx + 2t\int_0^\infty \frac{\tan^{-1}(x/t)}{(t^2+x^2)(e^{2\pi x}-1)}dx$$

$$= \frac{1}{2t}\log t - \frac{1}{2}t\left(2\frac{\log t}{t} - \frac{2}{t}\right) - \frac{1}{2}(\log^2 t - 2\log t + 2) - \frac{1}{2}\frac{\partial}{\partial t}\varsigma''(0,t)$$

$$= \frac{1}{2t}\log t - \frac{1}{2}\log^2 t - \frac{1}{2}\frac{\partial}{\partial t}\varsigma''(0,t)$$

We then have

$$\frac{1}{2}\frac{\partial}{\partial t}\varsigma''(0,t) = \frac{1}{2t}\log t - \frac{1}{2}\log^2 t + \int_0^\infty \frac{x\log(t^2+x^2)}{(t^2+x^2)(e^{2\pi x}-1)}dx - 2t\int_0^\infty \frac{\tan^{-1}(x/t)}{(t^2+x^2)(e^{2\pi x}-1)}dx$$

and noting the similarity with (5.1) we obtain

$$\frac{1}{2}\frac{\partial}{\partial t}\varsigma''(0,t) = \gamma_1(t)$$

We recall Lerch's identity [11]

(5.4)  $\log\Gamma(t) = \varsigma'(0,t) + \frac{1}{2}\log(2\pi)$

and note that differentiation of this gives us a similar relationship for $\gamma_0(t)$

$$\psi(t) = -\gamma_0(t) = \frac{\partial}{\partial t}\varsigma'(0,t)$$

In fact these equations are specific cases of the general result [22]

$$\int_1^t \gamma_n(x)dx = \frac{(-1)^{n+1}}{n+1}\left[\varsigma^{(n+1)}(0,t) - \varsigma^{(n+1)}(0)\right]$$

whereupon differentiation gives us

$$\gamma_n(t) = \frac{(-1)^{n+1}}{n+1}\frac{\partial}{\partial t}\varsigma^{(n+1)}(0,t)$$

We recall (1.5) which may be written as

$$\gamma_n(u) = \frac{1}{2u}\log^n u - \frac{1}{n+1}\log^{n+1} u - 2\operatorname{Re}\int_0^\infty \frac{i\log^n(u-ix)}{(u-ix)(e^{2\pi x}-1)}dx$$

where integration results in

$$\int_1^t \gamma_n(x)dx = \frac{1}{2(n+1)}\log^{n+1} t + (-1)^n n!t\sum_{k=0}^{n+1}\frac{(-1)^k \log^k t}{k!} - \frac{2}{n+1}\operatorname{Re} i\int_0^\infty \frac{\log^{n+1}(t-ix) - \log^{n+1}(1-ix)}{e^{2\pi x}-1}dx$$

where we have used

$$\int \log^{n+1} u\, du = (-1)^{n+1}(n+1)!u\sum_{k=0}^{n+1}\frac{(-1)^k \log^k u}{k!}$$

We then have

$$\frac{(-1)^{n+1}}{n+1}\left[\varsigma^{(n+1)}(0,t) - \varsigma^{(n+1)}(0)\right]$$

$$= \frac{1}{2(n+1)}\log^{n+1} t + (-1)^n n!t\sum_{k=0}^{n+1}\frac{(-1)^k \log^k t}{k!} - \frac{2}{n+1}\operatorname{Re} i\int_0^\infty \frac{\log^{n+1}(t-ix) - \log^{n+1}(1-ix)}{e^{2\pi x}-1}dx$$

and we note that this may also be obtained by differentiating (3.2) $n+1$ times.

It may be noted from (5.1) that

$$\lim_{u\to\infty}\left[\gamma_1(u) + \frac{1}{2}\log^2 u\right] = 0$$

which may be compared with (2.20).

Ivić [32, p.41] reports that the coefficients $\gamma_p(u)$ for $u \in (0,1]$ may be expressed as

$$\gamma_p(u) = \lim_{N \to \infty} \left( \sum_{n=0}^{N} \frac{\log^p(n+u)}{n+u} - \frac{\log^{p+1}(N+u)}{p+1} \right)$$

and we therefore have

$$\lim_{u \to \infty} \left[ \gamma_1(u) + \frac{1}{2}\log^2 u \right] = \lim_{u \to \infty} \left[ \lim_{N \to \infty} \left( \sum_{n=0}^{N} \frac{\log(n+u)}{n+u} - \frac{1}{2}\log^2(N+u) \right) + \frac{1}{2}\log^2 u \right] = 0$$

**6. Other integrals involving the Barnes multiple gamma functions**

Integrating (2.4) gives us

(6.1) $\quad \int_0^t \log \Gamma(u) du = -\frac{3}{4}t^2 - \frac{1}{2}t(1-t)\log t + \frac{1}{2}t + \frac{1}{2}t\log(2\pi) + 2\int_0^t du \int_0^\infty \frac{\tan^{-1}(x/u)}{e^{2\pi x}-1} dx$

and using (2.17)

$$2\int_0^v du \int_0^\infty \frac{\tan^{-1}(x/u)}{e^{2\pi x}-1} dx = 2v\int_0^\infty \frac{\tan^{-1}(x/v)}{e^{2\pi x}-1} dx + \int_0^\infty \frac{x\log(v^2+x^2)}{e^{2\pi x}-1} dx - 2\int_0^\infty \frac{x\log x}{e^{2\pi x}-1} dx$$

we deduce that

(6.2) $\quad \int_0^t \log \Gamma(u) du = -\frac{3}{4}t^2 - \frac{1}{2}t(1-t)\log t + \frac{1}{2}t + \frac{1}{2}t\log(2\pi) + 2t\int_0^\infty \frac{\tan^{-1}(x/t)}{e^{2\pi x}-1} dx$

$$+ \int_0^\infty \frac{x\log(t^2+x^2)}{e^{2\pi x}-1} dx - 2\int_0^\infty \frac{x\log x}{e^{2\pi x}-1} dx$$

An evaluation of the following well-known integral is contained in [22]

(6.3) $\quad \int_0^\infty \frac{x\log x}{e^{2\pi x}-1} dx = \frac{1}{2}\varsigma'(-1)$

With regard to the left-hand side of (6.2) we note Alexeiewsky's theorem [46, p.32]

(6.4) $\quad \int_0^t \log \Gamma(u) du = \frac{1}{2}t(1-t) + \frac{1}{2}t\log(2\pi) - \log G(1+t) + t\log \Gamma(t)$

where $G(x)$ is the Barnes double gamma function $\Gamma_2(x) = 1/G(x)$ defined, inter alia, by the Weierstrass canonical product [46, p.25]

(6.5) $$G(1+x) = (2\pi)^{x/2} \exp\left[-\frac{1}{2}(\gamma x^2 + x^2 + x)\right] \prod_{k=1}^{\infty} \left\{\left(1+\frac{x}{k}\right)^k \exp\left(\frac{x^2}{2k} - x\right)\right\}$$

and we note that $G(1) = G(2) = 1$.

We then have from (6.2) and (6.4)

(6.6) $$\frac{1}{2}t(1-t) - \log G(1+t) + t \log \Gamma(t) =$$

$$-\frac{3}{4}t^2 - \frac{1}{2}t(1-t)\log t + \frac{1}{2}t + 2t\int_0^\infty \frac{\tan^{-1}(x/t)}{e^{2\pi x} - 1} dx + \int_0^\infty \frac{x \log(t^2 + x^2)}{e^{2\pi x} - 1} dx - \varsigma'(-1)$$

and, using Binet's second formula (2.4) for $\log \Gamma(t)$, this becomes

(6.7) $$\log G(1+t) = \frac{1}{2}t^2\left[\log t - \frac{3}{2}\right] + \frac{1}{2}t\log(2\pi) + \varsigma'(-1) - \int_0^\infty \frac{x \log(t^2 + x^2)}{e^{2\pi x} - 1} dx$$

This result was previously obtained by Adamchik [3] in 2004. With $t = 0$ we recover equation (6.3), and with $t = 1$ we obtain

(6.8) $$\int_0^\infty \frac{x \log(1 + x^2)}{e^{2\pi x} - 1} dx = \frac{1}{2}\log(2\pi) + \varsigma'(-1) - \frac{3}{4}$$

which is also derived in [3] and [12] (along with several generalised integrals) using Hermite's integral for the Hurwitz zeta function.

We may also use (2.19) to determine that

(6.9) $$\log G(1+t) = t \log \Gamma(t) + \frac{1}{4}t^2 - \frac{1}{2}t(t-1)\log t - \int_0^\infty \frac{1-e^{-ty}}{y^2}\left[\frac{1}{e^y - 1} - \frac{1}{y} + \frac{1}{2}\right] dy$$

We also note the functional relationship with the Hurwitz zeta function due to Gosper [28] and Vardi [47] (a different derivation of which is contained in [21])

(6.10) $$\log G(1+t) - t \log \Gamma(t) = \varsigma'(-1) - \varsigma'(-1, t)$$

and we may therefore write (6.6) as

(6.11) $$\varsigma'(-1,t) = \frac{1}{2}t(t-1)\log t - \frac{1}{4}t^2 + 2t\int_0^\infty \frac{\tan^{-1}(x/t)}{e^{2\pi x}-1}dx + \int_0^\infty \frac{x\log(t^2+x^2)}{e^{2\pi x}-1}dx$$

as previously noted by Adamchik [3].

It is interesting to note that the Gosper/Vardi functional relationship may also be derived by comparing (3.5) with (6.7).

$\square$

Adamchik [3] showed that for $\text{Re}(v) > 0$ (where I have corrected the sign immediately preceding the integral)

(6.11.1)
$$\log G(1+v) = v\log\Gamma(v) + \frac{1}{4}v^2 - \frac{1}{2}B_2(v)\log v - \frac{1}{12} + \varsigma'(-1) + \int_0^\infty \frac{e^{-vt}}{t^2}\left[\frac{1}{1-e^{-t}} - \frac{1}{t} - \frac{1}{2} - \frac{1}{12}t\right]dt$$

$$= v\log\Gamma(v) + \frac{1}{4}v^2 - \frac{1}{2}B_2(v)\log v - \frac{1}{12} + \varsigma'(-1) + \int_0^\infty \frac{e^{-vt}}{t^2}\left[\frac{1}{e^t-1} - \frac{1}{t} + \frac{1}{2} - \frac{1}{12}t\right]dt$$

where $G(v)$ is the Barnes double gamma function defined in (6.5) below.

Letting $v = 1$ we obtain

$$\int_0^\infty \frac{e^{-t}}{t^2}\left[\frac{1}{e^t-1} - \frac{1}{t} + \frac{1}{2} - \frac{1}{12}t\right]dt = -\varsigma'(-1) - \frac{1}{6}$$

and we see that

$$v\log\Gamma(v) - \log G(1+v) + \frac{1}{4}v^2 - \frac{1}{2}B_2(v)\log v - \frac{1}{4} = \int_0^\infty \frac{e^{-t}-e^{-vt}}{t^2}\left[\frac{1}{e^t-1} - \frac{1}{t} + \frac{1}{2} - \frac{1}{12}t\right]dt$$

Using the well-known Frullani integral

$$\log v = \int_0^\infty \frac{e^{-t}-e^{-vt}}{t}dt$$

this becomes

$$v\log\Gamma(v) - \log G(1+v) + \frac{1}{4}v^2 - \frac{1}{2}\left[B_2(v) - \frac{1}{6}\right]\log v - \frac{1}{4} = \int_0^\infty \frac{e^{-t}-e^{-vt}}{t^2}\left[\frac{1}{e^t-1} - \frac{1}{t} + \frac{1}{2}\right]dt$$

Combining (2.19) and (6.7) gives us

$$\frac{1}{2}\log v - \frac{1}{2}v^2 \log v + \frac{1}{4}v^2 + v\log \Gamma(v) - \log G(1+v)$$

$$= \int_0^\infty \frac{1-e^{-vt}}{t^2}\left[\frac{1}{e^t-1} - \frac{1}{t} + \frac{1}{2}\right]dt$$

which agrees with (6.9)

We have

$$\int_0^\infty \frac{e^{-t}-e^{-vt}}{t^2}\left[\frac{1}{e^t-1} - \frac{1}{t} + \frac{1}{2}\right]dt = \int_0^\infty \frac{e^{-t}-1}{t^2}\left[\frac{1}{e^t-1} - \frac{1}{t} + \frac{1}{2}\right]dt + \int_0^\infty \frac{1-e^{-vt}}{t^2}\left[\frac{1}{e^t-1} - \frac{1}{t} + \frac{1}{2}\right]dt$$

and therefore we obtain

$$v\log \Gamma(v) - \log G(1+v) + \frac{1}{4}v^2 - \frac{1}{2}\left[B_2(v) - \frac{1}{6}\right]\log v - \frac{1}{4}$$

$$= \int_0^\infty \frac{e^{-t}-1}{t^2}\left[\frac{1}{e^t-1} - \frac{1}{t} + \frac{1}{2}\right]dt + \frac{1}{2}v\log v - \frac{1}{2}v^2 \log v + \frac{1}{4}v^2 - \log G(1+v) + v\log \Gamma(v)$$

We then see that

$$\int_0^\infty \frac{e^{-t}-1}{t^2}\left[\frac{1}{e^t-1} - \frac{1}{t} + \frac{1}{2}\right]dt = -\frac{1}{4}$$

which was derived earlier in (2.20).

Differentiating (6.11.1) gives us

$$\frac{G'(1+v)}{G(1+v)} = \log \Gamma(v) + v\psi(v) + \frac{1}{2}v - \frac{1}{2v}B_2(v) - \frac{1}{2}(2v-1)\log v - \int_0^\infty \frac{e^{-vt}}{t}\left[\frac{1}{e^t-1} - \frac{1}{t} + \frac{1}{2} - \frac{1}{12}t\right]dt$$

$$= \log \Gamma(v) + v\psi(v) + \frac{1}{2}v - \frac{1}{2v}B_2(v) - \frac{1}{2}(2v-1)\log v - \int_0^\infty \frac{e^{-vt}}{t}\left[\frac{1}{e^t-1} - \frac{1}{t} + \frac{1}{2}\right]dt + \frac{1}{12}\int_0^\infty e^{-vt}dt$$

and using (2.12) this becomes

$$= \log \Gamma(v) + v\psi(v) + \frac{1}{2}v - \frac{1}{2v}B_2(v) - \frac{1}{2}(2v-1)\log v + \left(v - \frac{1}{2}\right)\log v - v + \frac{1}{2}\log(2\pi) - \log \Gamma(v) + \frac{1}{12v}$$

Equality is easily verified using (6.24).

$\square$

Integrating (6.11) results in

$$\int_0^v \varsigma'(-1,t)\,dt = \left(\frac{1}{6}v^3 - \frac{1}{4}v^2\right)\log v - \frac{10}{72}v^3 + \frac{1}{8}v^2 + 2\int_0^v t\,dt \int_0^\infty \frac{\tan^{-1}(x/t)}{e^{2\pi x}-1}\,dx + \int_0^v dt \int_0^\infty \frac{x\log(t^2+x^2)}{e^{2\pi x}-1}\,dx$$

We have using integration by parts

$$\int t \tan^{-1}(x/t)\,dt = \frac{1}{2}t^2 \tan^{-1}(x/t) - \frac{1}{2}x^2 \tan^{-1}(t/x) + \frac{1}{2}xt$$

and since from (2.6) for $v/x > 0$

$$\tan^{-1}(x/v) = \frac{\pi}{2} - \tan^{-1}(v/x)$$

we obtain

$$\int_0^v t \tan^{-1}(x/t)\,dt = \frac{1}{2}(v^2 + x^2)\tan^{-1}(x/v) + \frac{1}{2}xv - \frac{\pi}{4}x^2$$

Integration by parts also gives us

$$\int_0^v \log(t^2 + x^2)\,dt = 2x\tan^{-1}(v/x) + v\log(v^2 + x^2) - 2v$$

$$= \pi x - 2x\tan^{-1}(x/v) + v\log(v^2 + x^2) - 2v$$

Therefore we have

$$I = 2\int_0^v t\,dt \int_0^\infty \frac{\tan^{-1}(x/t)}{e^{2\pi x}-1}\,dx + \int_0^v dt \int_0^\infty \frac{x\log(t^2+x^2)}{e^{2\pi x}-1}\,dx$$

$$= -\int_0^\infty \frac{x^2 \tan^{-1}(x/v)}{e^{2\pi x}-1}\,dx + v^2 \int_0^\infty \frac{\tan^{-1}(x/v)}{e^{2\pi x}-1}\,dx + v\int_0^\infty \frac{x\log(v^2+x^2)}{e^{2\pi x}-1}\,dx$$

$$-v\int_0^\infty \frac{x}{e^{2\pi x}-1}\,dx + \frac{\pi}{2}\int_0^\infty \frac{x^2}{e^{2\pi x}-1}\,dx$$

and using (2.4), (6.7) and (6.19) this becomes

$$I = -\int_0^\infty \frac{x^2 \tan^{-1}(x/v)}{e^{2\pi x}-1}\,dx + \frac{1}{2}v^2 \log\Gamma(v) - \frac{1}{2}v^2\left(v-\frac{1}{2}\right)\log v + \frac{1}{2}v^3 - \frac{1}{4}v^2 \log(2\pi)$$

(6.12) $$+ \frac{1}{2}v^3\left[\log v - \frac{3}{2}\right] + v\left[\frac{1}{2}v\log(2\pi) + \varsigma'(-1)\right] - v\log G(1+v)$$

$$- \frac{1}{24}v + \frac{\varsigma(3)}{8\pi^2}$$

We now need to evaluate $\int_0^v \varsigma'(-1,t)\,dt$ and this is carried out below.

Since the Hurwitz zeta function is analytic in the whole complex plane except for $s \neq 1$, its partial derivatives commute in the region where the function is analytic: we therefore have

$$\frac{\partial}{\partial t}\frac{\partial}{\partial s}\varsigma(s,t) = \frac{\partial}{\partial s}\frac{\partial}{\partial t}\varsigma(s,t) = -\frac{\partial}{\partial s}[s\varsigma(s+1,t)]$$

$$= -\varsigma(s+1,t) - s\frac{\partial}{\partial s}\varsigma(s+1,t)$$

and upon integrating with respect to $t$ we see that

$$-s\int_0^v \varsigma'(s+1,t)\,dt = \int_0^v \frac{\partial}{\partial t}\frac{\partial}{\partial s}\varsigma(s,t)\,dt + \int_0^v \varsigma(s+1,t)\,dt$$

We therefore get

$$-s\int_0^v \varsigma'(s+1,t)\,dt = \varsigma'(s,v) - \varsigma'(s,0) + \int_0^v \varsigma(s+1,t)\,dt$$

and with $s = -n$ we have

$$n\int_0^v \varsigma'(1-n,t)\,du = \varsigma'(-n,v) - \varsigma'(-n,0) + \int_0^v \varsigma(1-n,t)\,dt$$

Then, using the well-known result (3.10)

$$\varsigma(1-n,v) = -\frac{B_n(v)}{n} \text{ for } n \geq 1$$

we obtain

(6.13) $$n\int_0^v \varsigma'(1-n,t)\,dt = \frac{B_{n+1} - B_{n+1}(v)}{n(n+1)} + \varsigma'(-n,v) - \varsigma'(-n,0)$$

This integral was originally derived by Adamchik [2] in a different manner in 1998. With $n = 2$ we obtain

(6.14) $$\int_0^v \varsigma'(-1,t)\,dt = -\frac{1}{12}B_3(v) + \frac{1}{2}\varsigma'(-2,v) + \frac{\varsigma(3)}{8\pi^2}$$

since $\varsigma'(-n,0) = \varsigma'(-n)$ and $\varsigma'(-2) = -\frac{\varsigma(3)}{4\pi^2}$.

With $v = 1$ we note that

(6.15) $$\int_0^1 \varsigma'(1-n,t)\,dt = 0$$

Combining (6.12) and (6.14) gives us

(6.16) $$\int_0^\infty \frac{x^2 \tan^{-1}(x/v)}{e^{2\pi x} - 1}\,dx = \frac{1}{6}v^3 \log v - \frac{11}{36}v^3 + \frac{1}{2}v^2 \log \Gamma(v) + \frac{1}{4}v^2 \log(2\pi)$$

$$+v\varsigma'(-1) - \frac{1}{2}\varsigma'(-2,v) - v\log G(1+v)$$

With $v = 1$ we obtain

(6.17) $$\int_0^\infty \frac{x^2 \tan^{-1} x}{e^{2\pi x} - 1}\,dx = -\frac{11}{36} + \frac{1}{4}\log(2\pi) + \varsigma'(-1) + \frac{\varsigma(3)}{8\pi^2}$$

This concurs with the more general results previously obtained by Adamchik [3] and Boros et al. [12] using Hermite's integral (19).

We note in passing that the double gamma function may be represented by the logarithmic expansion [21]

$$\log G(1+t) - t \log \Gamma(t) = \varsigma'(-1) - \frac{1}{2}\varsigma(-1,t) - \frac{1}{2}\sum_{n=0}^{\infty} \frac{1}{n+1} \sum_{k=0}^{n} \binom{n}{k} (-1)^k (t+k)^2 \log(t+k)$$

We now integrate (6.7) to obtain

$$\int_0^u \log G(1+t)\,dt = \frac{1}{3}u^3 \left[\log u - \frac{3}{4}\right] + \frac{1}{4}u^2 \log(2\pi) + \varsigma'(-1)u - \int_0^u dt \int_0^\infty \frac{x \log(t^2 + x^2)}{e^{2\pi x} - 1}\,dx$$

Reversing the order of integration

$$\int_0^u dt \int_0^\infty \frac{x \log(t^2 + x^2)}{e^{2\pi x} - 1}\,dx = \int_0^\infty \frac{x}{e^{2\pi x} - 1}\,dx \int_0^u \log(t^2 + x^2)\,dt$$

and using integration by parts we have

$$\int_0^u \log(t^2 + x^2)\,dt = u \log(u^2 + x^2) + 2x \tan^{-1}(u/x) - 2u$$

This results in

(6.18) $$\int_0^u \log G(1+t)\,dt = \frac{1}{3}u^3 \left[\log u - \frac{3}{4}\right] + \left[\frac{1}{4}u^2 \log(2\pi) + \varsigma'(-1)u\right]$$

$$-u \int_0^\infty \frac{x \log(u^2 + x^2)}{e^{2\pi x} - 1}\,dx - 2 \int_0^\infty \frac{x^2 \tan^{-1}(u/x)}{e^{2\pi x} - 1}\,dx - 2u \int_0^\infty \frac{x}{e^{2\pi x} - 1}\,dx$$

With regard to the last integral, we have the well-known representation for the polylogarithm function $Li_s(z)$ [49, p.280] (a derivation is also contained in [22])

(6.19) $$Li_s(z) = \frac{z}{\Gamma(s)} \int_0^\infty \frac{u^{s-1}}{e^u - z}\,du$$

where $Li_s(z) = \sum_{n=1}^{\infty} \frac{z^n}{n^s}$ and $Li_s(1) = \varsigma(s)$.

Letting $z = 1$, $u = 2\pi x$ and $s = 2n$ we obtain

$$\frac{\varsigma(2n)\Gamma(2n)}{(2\pi)^{2n}} = \int_0^\infty \frac{x^{2n-1}}{e^{2\pi x} - 1} dx$$

Using Euler's identity for $\varsigma(2n)$ we therefore have an integral expression for the Bernoulli numbers

(6.20) $$B_{2n} = 4n(-1)^{n+1} \int_0^\infty \frac{x^{2n-1}}{e^{2\pi x} - 1} dx$$

and this corrects the typographical error in [13, p.223]. This formula was also derived by Ramanujan [10, Part II, p.220]. We therefore have the required particular case

(6.21) $$\int_0^\infty \frac{x}{e^{2\pi x} - 1} dx = \frac{B_2}{4} = \frac{1}{24}$$

We easily see from (6.19) that

$$\frac{\varsigma(s)\Gamma(s)}{(2\pi)^s} = \int_0^\infty \frac{x^{s-1}}{e^{2\pi x} - 1} dx$$

We note Legendre's relation (2.10)

$$2\int_0^\infty \frac{\sin(xy)}{e^{2\pi x} - 1} dx = \frac{1}{e^y - 1} - \frac{1}{y} + \frac{1}{2} = \frac{1}{2}\coth\frac{y}{2} - \frac{1}{y}$$

and, expressing $\sin(xy)$ as a power series, the left-hand side may be represented by

$$2\int_0^\infty \frac{\sin(xy)}{e^{2\pi x} - 1} dx = 2\sum_{n=1}^\infty \frac{(-1)^{n-1} y^{2n-1}}{(2n-1)!} \int_0^\infty \frac{x^{2n-1}}{e^{2\pi x} - 1} dx$$

The right-hand side may be expressed in terms of the Bernoulli numbers

$$\frac{y}{e^y - 1} = \sum_{n=0}^\infty B_n \frac{y^n}{n!} \qquad , (|y| < 2\pi)$$

Since $B_0 = 1$, $B_1 = -\frac{1}{2}$ and $B_{2n+1} = 0$ for all $n \geq 1$ we see that

$$\frac{1}{e^y-1} - \frac{1}{y} + \frac{1}{2} = \sum_{n=2}^{\infty} B_n \frac{y^{n-1}}{n!} = \sum_{n=1}^{\infty} B_{2n} \frac{y^{2n-1}}{(2n)!}$$

and equating coefficients of $y$ we deduce another derivation of (6.20).

We now substitute (6.7) and (6.21) in (6.18) to obtain

$$(6.22) \quad \int_0^u \log G(1+t)dt = \frac{1}{3}u^3\left[\log u - \frac{3}{4}\right] + \left[\frac{1}{4}u^2 \log(2\pi) + \varsigma'(-1)u\right]$$

$$+u\log G(1+u) - \frac{1}{2}u^3\left[\log u - \frac{3}{2}\right] + u\left[\frac{1}{2}u\log(2\pi) - \varsigma'(-1)\right] - 2\int_0^{\infty} \frac{x^2 \tan^{-1}(u/x)}{e^{2\pi x}-1}dx - \frac{1}{12}u$$

The following integral of the Barnes double gamma function is known from [46, p.207] (see also [21])

$$(6.23) \quad \int_0^u \log G(1+t)dt = \left[2\varsigma'(-1) - \frac{1}{12}\right]u + \frac{1}{4}u^2 \log(2\pi) - \frac{1}{6}u^3$$

$$+(u-1)\log G(1+u) + \log G(u) - 2\log \Gamma_3(1+u) + 2\log \Gamma_3(u)$$

where $\Gamma_3(u)$ is the triple gamma function defined in [46, p.42]

$$\Gamma_3(1+x) = \exp(c_1 x + c_2 x^2 + c_3 x^3)\prod_{k=1}^{\infty}\left\{\left(1+\frac{x}{k}\right)^{-\frac{1}{2}k(k+1)} \exp\left[\left(1+\frac{1}{k}\right)\left(\frac{1}{2}kx - \frac{1}{4}x^2 + \frac{1}{6k}x^3\right)\right]\right\}$$

where the constants are

$$c_1 = \frac{3}{8} - \frac{1}{4}\log(2\pi) - \log A, \quad c_2 = \frac{1}{4}\left[\gamma + \log(2\pi) + \frac{1}{2}\right]$$

$$c_3 = -\frac{1}{6}\left[\gamma + \varsigma(2) + \frac{3}{2}\right] \qquad \log A = \frac{1}{12} - \varsigma'(-1)$$

Hence, combining (6.22) and (6.23) will also enable us to evaluate the integral $\int_0^{\infty} \frac{x^2 \tan^{-1}(u/x)}{e^{2\pi x}-1}dx$ in terms of the multiple gamma functions.

We now multiply (6.7) by $t$ and integrate to obtain

$$\int_0^u t \log G(1+t)dt = \frac{1}{8}u^4\left[\log u - \frac{1}{4}u^4\right] - \frac{3}{16}u^4 - \left[\frac{1}{6}u^3 \log(2\pi) - \frac{1}{2}\varsigma'(-1)u^2\right]$$

$$-\int_0^u t\, dt \int_0^\infty \frac{x\log(t^2+x^2)}{e^{2\pi x}-1}dx$$

Reversing the order of integration

$$\int_0^u t\, dt \int_0^\infty \frac{x\log(t^2+x^2)}{e^{2\pi x}-1}dx = \int_0^\infty \frac{x}{e^{2\pi x}-1}dx \int_0^u t\log(t^2+x^2)dt$$

and we have

$$\int_0^u t\log(t^2+x^2)dt = \frac{1}{2}(u^2+x^2)\log(u^2+x^2) - \frac{1}{2}u^2 - x^2\log x$$

We then have

$$\int_0^u t \log G(1+t)dt = \frac{1}{8}u^4\left[\log u - \frac{1}{4}u^4\right] - \frac{3}{16}u^4 - \left[\frac{1}{6}u^3 \log(2\pi) - \frac{1}{2}\varsigma'(-1)u^2\right]$$

$$-\frac{1}{2}u^2\int_0^\infty \frac{x\log(u^2+x^2)}{e^{2\pi x}-1}dx - \frac{1}{2}\int_0^\infty \frac{x^3\log(u^2+x^2)}{e^{2\pi x}-1}dx$$

$$+\frac{1}{2}u^2\int_0^\infty \frac{x}{e^{2\pi x}-1}dx + \int_0^\infty \frac{x^3 \log x}{e^{2\pi x}-1}dx$$

We note from [46, p.209] that

$$\int_0^u t\log G(a+t)dt = (2-a)\left[-\frac{1}{4}+\frac{1}{2}(a-1)\log(2\pi) - 2\log A - \frac{1}{2}a^2 + a\right]u$$

$$+\frac{1}{2}\left[\frac{7}{4}+\frac{1}{2}\log(2\pi) - 2\log A + \frac{1}{2}a^2 - 2a\right]u^2 + \frac{1}{6}[\log(2\pi)-a]u^3$$

$$-\frac{1}{8}u^4 + [u^2 - a^2 + 4a - 4]\log G(u+a) + 2(a-2-u)\log \Gamma_3(u+a)$$

$$+(a-2)^2 \log G(a) + 2(2-a)\log \Gamma_3(a) + 2\int_0^u \log \Gamma_3(a+t)dt$$

and thus we get

$$\int_0^u t \log G(1+t)dt = \left[-\frac{1}{4} - 2\log A + \frac{1}{2}\right]u$$

$$+ \frac{1}{2}\left[\frac{7}{4} + \frac{1}{2}\log(2\pi) - 2\log A - \frac{3}{2}\right]u^2 + \frac{1}{6}[\log(2\pi)-1]u^3$$

$$-\frac{1}{8}u^4 + [u^2-1]\log G(1+u) - (1+u)\log \Gamma_3(1+u) + 2\int_0^u \log \Gamma_3(1+t)dt$$

Accordingly we may represent $\int_0^u \log \Gamma_3(1+t)dt$ in terms of $\int_0^\infty \frac{x^3 \log(u^2+x^2)}{e^{2\pi x}-1}dx$.

□

Alternatively, we now multiply (2.13) by $u$ and integrate to obtain

$$\int_0^v u\psi(u)\,du = -\frac{1}{2}v + \frac{1}{2}v^2 \log v - \frac{1}{4}v^2 - \int_0^\infty \frac{x\log(v^2+x^2)}{e^{2\pi x}-1}dx + 2\int_0^\infty \frac{x\log x}{e^{2\pi x}-1}dx$$

Integration by parts gives us

$$\int_0^v u\psi(u)\,du = v\log \Gamma(v) - \int_0^v \log \Gamma(u)\,du$$

Using (6.4) we obtain

$$\int_0^v u\psi(u)\,du = -\frac{1}{2}v(1-v) - \frac{1}{2}v\log(2\pi) + \log G(1+v)$$

and hence we have

$$-\frac{1}{2}v(1-v) - \frac{1}{2}v\log(2\pi) + \log G(1+v) = -\frac{1}{2}v + \frac{1}{2}v^2 \log v - \frac{1}{4}v^2 - \int_0^\infty \frac{x\log(v^2+x^2)}{e^{2\pi x}-1}dx + 2\int_0^\infty \frac{x\log x}{e^{2\pi x}-1}dx$$

Then using (6.3) this becomes (6.7).

We may also multiply (2.13) by $u^2$ and integrate to obtain

$$\int_0^v u^2 \psi(u)\,du = -\frac{1}{4}v^2 + \frac{1}{3}v^3 \log v - \frac{1}{9}v^3 - 2\int_0^v\int_0^\infty \left(1 - \frac{x^2}{u^2+x^2}\right)\frac{x}{e^{2\pi x}-1}\,dx\,du$$

$$\int_0^v\int_0^\infty \left(1 - \frac{x^2}{u^2+x^2}\right)\frac{x}{e^{2\pi x}-1}\,dx\,du = v\int_0^\infty \frac{x}{e^{2\pi x}-1}\,dx - \int_0^\infty \frac{x^2 \tan^{-1}(v/x)}{e^{2\pi x}-1}\,dx$$

$$= v\int_0^\infty \frac{x}{e^{2\pi x}-1}\,dx + \int_0^\infty \frac{x^2 \tan^{-1}(x/v)}{e^{2\pi x}-1}\,dx - \frac{\pi}{2}\int_0^\infty \frac{x^2}{e^{2\pi x}-1}\,dx$$

$$= \int_0^\infty \frac{x^2 \tan^{-1}(x/v)}{e^{2\pi x}-1}\,dx + \frac{1}{24}v - \frac{\varsigma(3)}{8\pi^2}$$

Since

$$\int_0^v u^2 \psi(u)\,du = v^2 \log \Gamma(v) - 2\int_0^v u \log \Gamma(u)\,du$$

this also enables us to evaluate $\int_0^\infty \frac{x^2 \tan^{-1}(x/v)}{e^{2\pi x}-1}\,dx$.

□

Differentiating (6.7) gives us

$$\frac{G'(1+t)}{G(1+t)} = 2t + t\log t + \frac{1}{2}\log(2\pi) - 2t\int_0^\infty \frac{x}{(t^2+x^2)(e^{2\pi x}-1)}\,dx$$

and using (2.13) we see that

(6.24) $$\frac{G'(1+t)}{G(1+t)} = \frac{1}{2}\log(2\pi) + \frac{1}{2} - t + t\psi(t)$$

as reported in [49, p.264].

**7. Ramanujan's formula**

It was noted in [22] that

(7.1) $$(-1)^r \varsigma^{(r)}(s,u) = \frac{r!}{(s-1)^{r+1}} + \sum_{n=0}^\infty \frac{(-1)^n}{n!}(s-1)^n \gamma_{n+r}(u)$$

When $r = 1$ and $s = 0$, and making use of Lerch's identity [11],

$$\log \Gamma(u) = \varsigma'(0, u) + \frac{1}{2}\log(2\pi)$$

this reverts to (2.2).

In 1995, Choudhury [17] mentioned that Ramanujan had determined that for $r \geq 1$ and $\operatorname{Re}(s) > 1$

(7.2) $\qquad (-1)^r \varsigma^{(r)}(s) = \sum_{n=1}^{\infty} \frac{\log^r n}{n^s} = \frac{r!}{(s-1)^{r+1}} + \sum_{n=0}^{\infty} \frac{(-1)^n}{n!}(s-1)^n \gamma_{n+r}$

(see also Ramanujan's Notebooks [10, Part I, p.224]). This is a particular case of (7.1) where $u = 1$ and it should be noted that (7.1) is in fact valid for all $s \neq 1$.

As before, substituting Coffey's representation (1.5) in (7.1) we obtain

$$\sum_{n=0}^{\infty} \frac{(-1)^n}{n!}(s-1)^n \gamma_{n+r}(u) = \frac{\log^r u}{2u} \sum_{n=0}^{\infty} \frac{(-1)^n (s-1)^n}{n!} \log^n u - \sum_{n=0}^{\infty} \frac{(-1)^n (s-1)^n}{n!} \frac{1}{n+r+1} \log^{n+r+1} u$$

$$-2\operatorname{Re}\int_0^{\infty} \frac{i(u-ix)\log^r(u-ix) \sum_{n=0}^{\infty} \frac{(-1)^n (s-1)^n}{n!}\log^n(u-ix)}{(u^2 + x^2)(e^{2\pi x} - 1)} dx$$

For convenience we define $I(r, s)$ as

$$I(r, s) = -2\operatorname{Re}\int_0^{\infty} \frac{i(u-ix)\log^r(u-ix) \sum_{n=0}^{\infty} \frac{(-1)^n (s-1)^n}{n!}\log^n(u-ix)}{(u^2 + x^2)(e^{2\pi x} - 1)} dx$$

and, using the same modus operandi as before, the integrand may be simplified to

$$\sum_{n=0}^{\infty} \frac{(-1)^n (s-1)^n}{n!} \log^n(u-ix) = \exp[-(s-1)\log(u-ix)] = \exp[\log \frac{1}{(u-ix)^{s-1}}] = \frac{1}{(u-ix)^{s-1}}$$

We thus obtain

$$I(r, s) = -2\operatorname{Re}\int_0^{\infty} \frac{i(u+ix)^s \log^r(u-ix)}{(u^2 + x^2)^s (e^{2\pi x} - 1)} dx$$

Letting $s = 1$ we obtain equation (1.5).

With $s = 0$ we have

$$I(r,0) = -2\operatorname{Re}\int_0^\infty \frac{i\log^r(u-ix)}{e^{2\pi x}-1}dx$$

and with $r = 1$ this becomes

$$I(1,0) = -2\operatorname{Re}\int_0^\infty \frac{i\log(u-ix)}{e^{2\pi x}-1}dx$$

Therefore we have

$$I(1,0) = -2\int_0^\infty \frac{\tan^{-1}(x/u)}{e^{2\pi x}-1}dx$$

When $s = -1$ we obtain

$$I(r,-1) = -2\operatorname{Re}\int_0^\infty \frac{i(u^2+x^2)\log^r(u-ix)}{(u+ix)(e^{2\pi x}-1)}dx = -2\operatorname{Re}\int_0^\infty \frac{i(u-ix)\log^r(u-ix)}{e^{2\pi x}-1}dx$$

Letting $r = 0$ in the above equation results in

$$I(0,-1) = -2\operatorname{Re}\int_0^\infty \frac{i(u-ix)}{e^{2\pi x}-1}dx = 2\int_0^\infty \frac{x}{e^{2\pi x}-1}dx = \frac{1}{12}$$

where we have used (6.21).

When $r = 1$ we obtain

$$I(1,-1) = -\int_0^\infty \frac{x\log(u^2+x^2)+2u\tan^{-1}(u/x)}{e^{2\pi x}-1}dx$$

and reference to (6.11)

$$\varsigma'(-1,t) = \frac{1}{2}t(t-1)\log t - \frac{1}{4}t^2 + 2t\int_0^\infty \frac{\tan^{-1}(x/t)}{e^{2\pi x}-1}dx + \int_0^\infty \frac{x\log(t^2+x^2)}{e^{2\pi x}-1}dx$$

then tells us that

$$I(1,-1) = \frac{1}{2}t(t-1)\log t - \frac{1}{4}t^2 - \varsigma'(-1,t)$$

## 8. Some connections with sine and cosine integrals

The following was posed as a question in Whittaker & Watson [49, p.261]: Prove that for all values of $u$ except negative real values

$$(8.1) \qquad \log \Gamma(u) = \left(u - \frac{1}{2}\right)\log u - u + \frac{1}{2}\log(2\pi) + \frac{1}{\pi}\sum_{n=1}^{\infty}\int_0^{\infty}\frac{\sin(2n\pi x)}{n(x+u)}dx$$

and this result was attributed by Stieltjes to Bourguet. Equation (8.1) may also be derived using the Euler-Maclaurin summation formula (see, in particular, Knopp's book [35, p.530]).

The sine and cosine integrals, $Si(x)$ and $Ci(x)$, are respectively defined by

$$Si(x) = \int_0^x \frac{\sin t}{t}dt$$

$$Ci(x) = -\int_x^{\infty} \frac{\cos t}{t}dt$$

We have the well-known integral from Fourier series analysis

$$\frac{\pi}{2} = \int_0^{\infty}\frac{\sin t}{t}dt$$

which, as in [9, p.273], may be formally obtained by letting $u \to 0$ in (2.7).

Therefore defining

$$si(x) = Si(x) - \frac{\pi}{2}$$

we have

$$si(x) = \int_0^x \frac{\sin t}{t}dt - \int_0^{\infty}\frac{\sin t}{t}dt = -\int_x^{\infty}\frac{\sin t}{t}dt$$

By differentiation we easily see that

$$\frac{d}{dx}\left(\cos(2n\pi u)Si[2n\pi(x+u)] - \sin(2n\pi u)Ci[2n\pi(x+u)]\right) = \frac{\sin(2n\pi x)}{x+u}$$

and we therefore have

$$\int_0^M \frac{\sin(2n\pi x)}{x+u} dx = \left(\cos(2n\pi u) Si[2n\pi(x+u)] - \sin(2n\pi u) Ci[2n\pi(x+u)]\right)\Big|_0^M$$

$$= \cos(2n\pi u)\{Si[2n\pi(M+u)] - Si[2n\pi u]\} - \sin(2n\pi u)\{Ci[2n\pi(M+u)] - Ci[2n\pi u]\}$$

From the above definitions we see that

$$\lim_{M\to\infty} Si[2n\pi(M+a)] = \frac{\pi}{2}$$

and

$$\lim_{M\to\infty} Ci[2n\pi(M+a)] = 0$$

Hence we obtain as $M \to \infty$

$$\int_0^\infty \frac{\sin(2n\pi x)}{x+u} dx = \cos(2n\pi u)\left\{\frac{\pi}{2} - Si(2n\pi u)\right\} + \sin(2n\pi u) Ci(2n\pi u)$$

$$= -\cos(2n\pi u) si(2n\pi u) + \sin(2n\pi u) Ci(2n\pi u)$$

Therefore from Bourguet's formula we have (as also shown in [24])

(8.2)
$$\log \Gamma(u) = \left(u - \frac{1}{2}\right)\log u - u + \frac{1}{2}\log(2\pi) + \frac{1}{\pi}\sum_{n=1}^\infty \frac{1}{n}[\sin(2n\pi u) Ci(2n\pi u) - \cos(2n\pi u) si(2n\pi u)]$$

which was also reported by Nörlund [38, p.114].

Comparing (2.4) with (8.1) we note that

(8.3) $$2\int_0^\infty \frac{\tan^{-1}(x/u)}{e^{2\pi x} - 1} dx = \frac{1}{\pi}\sum_{n=1}^\infty \int_0^\infty \frac{\sin(2n\pi x)}{n(x+u)} dx$$

We may integrate (8.2) using Alexeiewsky's theorem [46, p.32]

$$\int_0^x \log \Gamma(u) du = \frac{x(1-x)}{2} + \frac{x}{2}\log(2\pi) - \log G(x+1) + x\log \Gamma(x)$$

For the integral of the right-hand side of (8.2) we have (as shown in equation (6.117fi) in [24])

$$\int_0^x \log \Gamma(u)\,du = \frac{1}{2}x\log(2\pi) + \frac{1}{4}x\left[2 - x + 2(x-1)\log x\right] - \frac{1}{2}x^2 + \frac{1}{12}\log x + \frac{1}{12} - \varsigma'(-1)$$

$$- \frac{1}{2\pi^2}\sum_{n=1}^{\infty}\frac{1}{n^2}\left[\cos(2n\pi x)Ci(2n\pi x) + \sin(2n\pi x)si(2n\pi x)\right]$$

We then determine that

(8.4) $\quad x\log\Gamma(x) - \log G(x+1) = \frac{1}{4}x\left[2(x-1)\log x\right] - \frac{1}{4}x^2 + \frac{1}{12}\log x + \frac{1}{12} - \varsigma'(-1)$

$$- \frac{1}{2\pi^2}\sum_{n=1}^{\infty}\frac{1}{n^2}\left[\cos(2n\pi x)Ci(2n\pi x) + \sin(2n\pi x)si(2n\pi x)\right]$$

We now recall (6.10)

$$\log G(1+x) - x\log\Gamma(x) = \varsigma'(-1) - \varsigma'(-1,x)$$

and deduce that

(8.5)
$$\varsigma'(-1,x) = -\varsigma(-1,x)\log x - \frac{1}{4}x^2 + \frac{1}{12} - \frac{1}{2\pi^2}\sum_{n=1}^{\infty}\frac{1}{n^2}\left[\cos(2n\pi x)Ci(2n\pi x) + \sin(2n\pi x)si(2n\pi x)\right]$$

where $\varsigma'(s,u) = \frac{\partial}{\partial s}\varsigma(s,u)$. The above formula was reported by Elizalde [27] in 1985.

Further exploratory work on the intimate relationship between the sine and cosine integrals and the Riemann zeta function is contained in [24].

It was shown in [21] that

$$\log\Gamma(u) = \sum_{n=0}^{\infty}\frac{1}{n+1}\sum_{k=0}^{n}\binom{n}{k}(-1)^k(u+k)\log(u+k) + \frac{1}{2} - u + \frac{1}{2}\log(2\pi)$$

and comparing this with (2.4) we see that

$$2\int_0^{\infty}\frac{\tan^{-1}(x/u)}{e^{2\pi x}-1}dx = \sum_{n=0}^{\infty}\frac{1}{n+1}\sum_{k=0}^{n}\binom{n}{k}(-1)^k(u+k)\log(u+k) + \frac{1}{2} - \left(u - \frac{1}{2}\right)\log u$$

Letting $x = uv$ in (2.13) results in

$$\psi(u) = -\frac{1}{2u} + \log u - 2\int_0^\infty \frac{v}{(1+v^2)(e^{2\pi uv}-1)} dv$$

We now let $u$ be a positive integer

$$\psi(n) = -\frac{1}{2n} + \log n - 2\int_0^\infty \frac{v}{(1+v^2)(e^{2\pi nv}-1)} dv$$

and make the summation

$$\sum_{n=1}^\infty \frac{\psi(n)}{n^s} = -\frac{1}{2}\sum_{n=1}^\infty \frac{1}{n^{s+1}} + \sum_{n=1}^\infty \frac{\log n}{n^s} - 2\sum_{n=1}^\infty \frac{1}{n^s}\int_0^\infty \frac{v}{(1+v^2)(e^{2\pi nv}-1)} dv$$

$$= -\frac{1}{2}\varsigma(s+1) - \varsigma'(s) - 2\sum_{n=1}^\infty \frac{1}{n^s}\int_0^\infty \frac{v}{(1+v^2)(e^{2\pi nv}-1)} dv$$

Assuming that it is valid to interchange the order of summation and integration we have

$$\sum_{n=1}^\infty \frac{\psi(n)}{n^s} = -\frac{1}{2}\varsigma(s+1) - \varsigma'(s) - 2\int_0^\infty \frac{v}{1+v^2}\sum_{n=1}^\infty \frac{1}{n^s(e^{2\pi nv}-1)} dv$$

Ogreid and Osland [39] report that

(8.6) $$\sum_{n=1}^\infty \frac{\psi(n)}{n^2} = \varsigma(3) - \gamma\varsigma(2)$$

and we therefore end up with

$$\gamma\varsigma(2) - \frac{3}{2}\varsigma(3) - \varsigma'(2) = 2\int_0^\infty \frac{v}{1+v^2}\sum_{n=1}^\infty \frac{1}{n^2(e^{2\pi nv}-1)} dv$$

Plouffe (see [40] and [41]) has reported that

$$\varsigma(3) = \frac{7\pi^3}{180} - 2\sum_{n=1}^\infty \frac{1}{n^3(e^{2\pi n}-1)}$$

Bradley [14] also deals with similar series in his paper entitled "Ramanujan's formula for the logarithmic derivative of the gamma function".

Plouffe's representation lead me to suspect that perhaps $\sum_{n=1}^{\infty} \frac{1}{n^s(e^{2\pi nv}-1)}$ was related to the polylogarithm function $Li_s(v)$ and indeed, this was made clear in the recent paper by Vepstas [48] who noted that

$$\sum_{n=1}^{\infty} \frac{1}{n^s(e^{2\pi nv}-1)} = \sum_{n=1}^{\infty}\sum_{m=0}^{\infty} \frac{e^{-2\pi nv(m+1)}}{n^s} = \sum_{m=1}^{\infty} Li_s(e^{-2\pi mv})$$

We then have

$$\int_0^{\infty} \frac{v}{1+v^2} \sum_{n=1}^{\infty} \frac{1}{n^2(e^{2\pi nv}-1)} dv = \int_0^{\infty} \frac{v}{1+v^2} \sum_{m=1}^{\infty} Li_2(e^{-2\pi mv}) dv$$

$$= \sum_{m=1}^{\infty}\sum_{n=1}^{\infty} \frac{1}{n^2} \int_0^{\infty} \frac{ve^{-2\pi nmv}}{1+v^2} dv$$

We have from G&R [29, p.338]

(8.7) $$\int_0^{\infty} \frac{ve^{-av}}{b^2+v^2} dv = -[Ci(ba)\cos(ba) + si(ba)\sin(ba)]$$

and, in particular, we obtain

$$\int_0^{\infty} \frac{ve^{-2\pi nmv}}{1+v^2} dv = -Ci(2\pi nm)\cos(2\pi nm) = -Ci(2\pi nm)$$

Therefore we have

$$\int_0^{\infty} \frac{v}{1+v^2} \sum_{n=1}^{\infty} \frac{1}{n^2(e^{2\pi nv}-1)} dv = -\sum_{m=1}^{\infty}\sum_{n=1}^{\infty} \frac{1}{n^2} Ci(2\pi nm)$$

and thus

$$\gamma\varsigma(2) - \frac{3}{2}\varsigma(3) - \varsigma'(2) = -2\sum_{m=1}^{\infty}\sum_{n=1}^{\infty} \frac{1}{n^2} Ci(2\pi nm)$$

$$= -2\sum_{n=1}^{\infty}\frac{1}{n^2}\sum_{m=1}^{\infty}Ci(2\pi nm)$$

We now use an identity which appears in Nörlund's book [38, p.108] (see also equation (6.117q) of [24])

$$\psi(a) = \log a - \frac{1}{2a} + 2\sum_{n=1}^{\infty}[\cos(2n\pi a)Ci(2n\pi a) + \sin(2n\pi a)si(2n\pi a)]$$

I subsequently noted that the above identity may be derived directly from (2.13) and (8.7) (8.2), or alternatively by differentiating (8.2).

Letting $a = n$ results in

$$\psi(n) = \log n - \frac{1}{2n} + 2\sum_{m=1}^{\infty}Ci(2\pi nm)$$

and we therefore have

$$-2\sum_{n=1}^{\infty}\frac{1}{n^2}\sum_{m=1}^{\infty}Ci(2\pi nm) = \sum_{n=1}^{\infty}\frac{\log n}{n^2} - \sum_{n=1}^{\infty}\frac{\psi(n)}{n^2} - \frac{1}{2}\varsigma(3)$$

$$= \gamma\varsigma(2) - \frac{3}{2}\varsigma(3) - \varsigma'(2)$$

which simply leads us back to square one!

As a matter of interest, Abramowitz and Stegun [1, p.232] define auxiliary functions

$$f(x) = -\cos x\, si(x) + \sin x\, Ci(x) = \int_0^{\infty}\frac{\sin y}{y+x}dy$$

$$g(x) = -\cos x\, Ci(x) - \sin x\, si(x) = \int_0^{\infty}\frac{\cos y}{y+x}dy$$

and report that for $\mathrm{Re}(x) > 0$

$$f(x) = \int_0^{\infty}\frac{e^{-xu}}{1+u^2}du$$

$$g(x) = \int_0^\infty \frac{u e^{-xu}}{1+u^2} du$$

The above results may be derived by considering the iterated integral

$$I = \int_0^\infty \int_0^\infty e^{-(a+y)x} \sin y \, dx \, dy$$

where integrating with respect to $x$ gives us

$$\int_0^\infty e^{-(a+y)x} dx = \frac{1}{a+y}$$

and thus we have

$$I = \int_0^\infty \frac{\sin y}{a+y} dy$$

Similarly, integrating with respect to $y$ gives us

$$\int_0^\infty e^{-(a+y)x} \sin y \, dy = \frac{e^{-ax}}{2i} \int_0^\infty e^{-yx}(e^{iy} - e^{-iy}) dy = \frac{e^{-ax}}{1+x^2}$$

Therefore we see that

$$\int_0^\infty \frac{e^{-ax}}{1+x^2} dx = \int_0^\infty \frac{\sin y}{a+y} dy$$

and the validity of the operation

$$\int_0^\infty dx \int_0^\infty e^{-(a+y)x} \sin y \, dy = \int_0^\infty dy \int_0^\infty e^{-(a+y)x} \sin y \, dx$$

is confirmed by [9, p.282]. The formula

$$\int_0^\infty \frac{x e^{-ax}}{1+x^2} dx = \int_0^\infty \frac{\cos y}{a+y} dy$$

may be derived in a similar fashion.

Letting $t = xy$ we see that

$$\int_0^\infty \frac{\sin(xy)}{x+a}\,dx = \int_0^\infty \frac{\sin t}{t+ay}\,dt$$

We may also make the summation for $\operatorname{Re}(p) > 1$

$$\sum_{n=1}^\infty \frac{1}{n^p}\int_0^\infty \frac{e^{-nxu}}{1+u^2}\,du = \sum_{n=1}^\infty \frac{1}{n^p}[-\cos(nx)si(nx) + \sin(nx)Ci(nx)]$$

and therefore we obtain

$$\int_0^\infty \frac{Li_p\left[e^{-xu}\right]}{1+u^2}\,du = \sum_{n=1}^\infty \frac{1}{n^p}[-\cos(nx)si(nx) + \sin(nx)Ci(nx)]$$

Similarly we have

$$\int_0^\infty \frac{uLi_p\left[e^{-xu}\right]}{1+u^2}\,du = -\sum_{n=1}^\infty \frac{1}{n^p}[\sin(nx)si(nx) + \cos(nx)Ci(nx)]$$

In the case $p = 2$, it may be noted from equation (6.117k) in [24] that the above integral may also be evaluated in terms of the gamma function and the Barnes double gamma function.

Differentiating (2.13) results in

$$\psi'(u) = \frac{1}{2u^2} + \frac{1}{u} + 4u\int_0^\infty \frac{x}{(u^2+x^2)^2(e^{2\pi x}-1)}\,dx$$

$$= \frac{1}{2u^2} + \frac{1}{u} + 4u\sum_{n=1}^\infty \int_0^\infty \frac{xe^{-2n\pi x}}{(u^2+x^2)^2}\,dx$$

and we note from G&R [29, p.338] that

$$\int_0^\infty \frac{xe^{-ax}}{(u^2+x^2)^2}\,dx = -\frac{a}{2u}[Ci(au)\sin(au) - si(au)\cos(au)]$$

Hence we see that (albeit the convergence of the series needs to be established)

$$\psi'(u) = \frac{1}{2u^2} + \frac{1}{u} - 4\pi\sum_{n=1}^\infty n[Ci(2n\pi u)\sin(2n\pi u) - si(2n\pi u)\cos(2n\pi u)]$$

Differentiating (2.13) gives us

$$\psi'(u) = \frac{1}{2u^2} + \frac{1}{u} + 4\pi u \int_0^\infty \frac{v e^{2\pi uv}}{(1+v^2)(e^{2\pi uv}-1)^2} dv$$

and hence we see that

$$\int_0^\infty \frac{v}{(u^2+v^2)^2(e^{2\pi v}-1)} dv = \pi \int_0^\infty \frac{v e^{2\pi uv}}{(1+v^2)(e^{2\pi uv}-1)^2} dv$$

With $u = 1$ we get

$$\int_0^\infty \frac{v}{(1+v^2)^2(e^{2\pi v}-1)} dv = \pi \int_0^\infty \frac{v e^{2\pi v}}{(1+v^2)(e^{2\pi v}-1)^2} dv$$

In passing, we note from G&R [29, p.428, eq. 3.749.2]

$$\int_0^\infty \frac{y \cot(\pi y)}{y^2 + x^2} dy = \frac{\pi}{e^{2\pi x} - 1}$$

## 9. Some double integrals for the Stieltjes constants

Kanemitsu et al. [38] showed in 2004 that

(9.1) $$S = \sum_{n=1}^\infty H_n \left( \log \frac{n+1}{n} - \frac{1}{n} \right) = -\frac{1}{2}\left[ \varsigma(2) + \gamma^2 - 2\gamma_1 \right]$$

and we now wish to represent this series as an integral. We have the well-known integral

$$\log t = \int_0^1 \frac{y^{t-1} - 1}{\log y} dy$$

which may be easily derived by integrating $\frac{1}{x} = \int_0^1 t^{x-1} dt$, and hence we have

$$\log \frac{n+1}{n} = \int_0^1 \frac{y^n - y^{n-1}}{\log y} dy = \int_0^1 \frac{y^{n-1}(y-1)}{\log y} dy$$

We also have the known integral expression for the harmonic numbers $H_n$ (see for example [21])

$$H_n = -n\int_0^1 (1-x)^{n-1} \log x\, dx$$

We therefore obtain

$$S = \sum_{n=1}^{\infty}\left[ n\int_0^1 (1-x)^{n-1}\log x\,dx \int_0^1 \frac{y^{n-1}(1-y)}{\log y}\,dy + \int_0^1 (1-x)^{n-1}\log x\,dx \right]$$

$$= \sum_{n=1}^{\infty}\left[ n\int_0^1 (1-x)^{n-1}\log x\,dx \int_0^1 \frac{y^{n-1}(1-y)}{\log y}\,dy + \int_0^1 (1-x)^{n-1}\log x\,dx \int_0^1 dy \right]$$

Using the derivative of the geometric series we have

$$\sum_{n=1}^{\infty} n\, z^{n-1} = \frac{1}{(1-z)^2}$$

and hence we obtain

$$S = \int_0^1\int_0^1 \left[ \frac{(1-y)\log x}{[1-(1-x)y]^2 \log y} - \frac{\log x}{x} \right] dx\,dy$$

We therefore get

(9.2) $$\int_0^1\int_0^1 \left[ \frac{(1-y)\log x}{[1-(1-x)y]^2 \log y} - \frac{\log x}{x} \right] dx\,dy = -\frac{1}{2}\left[\varsigma(2) + \gamma^2 - 2\gamma_1\right]$$

□

We may also write the double integral

$$\log^2 t = \int_0^1 \frac{x^{t-1}-1}{\log x}\,dx \int_0^1 \frac{y^{t-1}-1}{\log y}\,dy$$

$$= \int_0^1\int_0^1 \frac{(xy)^{t-1} - x^{t-1} - y^{t-1} + 1}{\log x \log y}\,dx\,dy$$

and we note from (2.1) that

$$\gamma_1(u) = -\frac{1}{2}\sum_{k=0}^{\infty} \frac{1}{k+1}\sum_{j=0}^{k}\binom{k}{j}(-1)^j \log^2(u+j)$$

We then have

$$\sum_{j=0}^{k}\binom{k}{j}(-1)^{j}\log^{2}(u+j)=\sum_{j=0}^{k}\binom{k}{j}(-1)^{j}\int_{0}^{1}\int_{0}^{1}\frac{(xy)^{u-1}(xy)^{j}-x^{u-1}x^{j}-y^{u-1}y^{j}+1}{\log x\log y}\,dx\,dy$$

$$=\int_{0}^{1}\int_{0}^{1}\frac{(xy)^{u-1}(1-xy)^{k}-x^{u-1}(1-x)^{k}-y^{u-1}(1-y)^{k}+\delta_{k,0}}{\log x\log y}\,dx\,dy$$

and we therefore obtain using the Maclaurin expansion for $\log(1-u)$

(9.3) $\quad \gamma_{1}(u)=\dfrac{1}{2}\displaystyle\int_{0}^{1}\int_{0}^{1}\left[\dfrac{(xy)^{u-1}\log(xy)}{(1-xy)\log x\log y}-\dfrac{x^{u-1}}{(1-x)\log y}-\dfrac{y^{u-1}}{(1-y)\log x}-\dfrac{1}{\log x\log y}\right]dx\,dy$

With $u=1$ we have

(9.4) $\quad \gamma_{1}=\dfrac{1}{2}\displaystyle\int_{0}^{1}\int_{0}^{1}\left[\dfrac{\log(xy)}{(1-xy)\log x\log y}-\dfrac{1}{(1-x)\log y}-\dfrac{1}{(1-y)\log x}-\dfrac{1}{\log x\log y}\right]dx\,dy$

and this complements the double integral representation previously found for Euler's constant $\gamma$ by Guillera and Sondow [30] in 2005

(9.5) $\quad \gamma=-\displaystyle\int_{0}^{1}\int_{0}^{1}\dfrac{1-x}{(1-xy)\log(xy)}\,dx\,dy$

and hence we have a particular value for the integral in (9.3).

We note, unlike the expression for $\gamma_{1}$, that (9.5) is not symmetrical in $x$ and $y$: symmetry may of course be easily restored by writing the double integral in the form

(9.6) $\quad \gamma=-\dfrac{1}{2}\displaystyle\int_{0}^{1}\int_{0}^{1}\dfrac{2-x-y}{(1-xy)\log(xy)}\,dx\,dy$

We note from equation (4.3.263) in [22] that

$$\gamma_{1}\left(\frac{1}{2}\right)=\gamma_{1}-\log^{2}2-2\gamma\log 2$$

Differentiating (9.3) gives us

$$\gamma_1'(u) = \frac{1}{2}\int_0^1\int_0^1 \left[\frac{(xy)^{u-1}\log^2(xy)}{(1-xy)\log x \log y} - \frac{x^{u-1}\log x}{(1-x)\log y} - \frac{y^{u-1}\log y}{(1-y)\log x}\right] dx\, dy$$

and we have

$$\gamma_1'(1) = \frac{1}{2}\int_0^1\int_0^1 \left[\frac{\log^2(xy)}{(1-xy)\log x \log y} - \frac{\log x}{(1-x)\log y} - \frac{\log y}{(1-y)\log x}\right] dx\, dy$$

We note from equation (4.3.244) in [22] that

$$\gamma_1'(1) = 2\pi^2 \varsigma'(-1) + \varsigma(2)(\gamma + \log 2\pi)$$

Using the integral

$$\log t = \int_0^1 \frac{x^{t-1}-1}{\log x} dx$$

and noting from (2.1) that

$$\gamma_0(u) = -\sum_{k=0}^\infty \frac{1}{k+1}\sum_{j=0}^k \binom{k}{j}(-1)^j \log(u+j)$$

we then have

$$\sum_{j=0}^k \binom{k}{j}(-1)^j \log(u+j) = \sum_{j=0}^k \binom{k}{j}(-1)^j \int_0^1 \frac{x^{u-1}x^j - 1}{\log x} dx$$

$$= \int_0^1 \frac{x^{u-1}(1-x)^k - \delta_{k,0}}{\log x} dx$$

Using the Maclaurin expansion for $\log(1-u)$ we therefore obtain the well-known integral

$$\gamma_0(u) = -\psi(u) = \int_0^1 \left[\frac{x^{u-1}}{1-x} + \frac{1}{\log x}\right] dx$$

☐

Nielsen [37, p.52] reports that

(9.4) $$[\psi(x)+\gamma]^2 = \psi'(x) - \varsigma(2) - 2\xi(x)$$

where

$$\xi(x) = \sum_{n=1}^{\infty} H_n \left( \frac{1}{x+n} - \frac{1}{n+1} \right)$$

Integration results in

$$\int_0^1 \xi(x)\, dx = \sum_{n=1}^{\infty} H_n \left( \log\frac{n+1}{n} - \frac{1}{n+1} \right)$$

$$= \sum_{n=1}^{\infty} H_n \left( \log\frac{n+1}{n} - \frac{1}{n} + \frac{1}{n} - \frac{1}{n+1} \right)$$

$$= \sum_{n=1}^{\infty} H_n \left( \log\frac{n+1}{n} - \frac{1}{n} \right) + \sum_{n=1}^{\infty} H_n \left( \frac{1}{n} - \frac{1}{n+1} \right)$$

and we easily see that

$$\sum_{n=1}^{\infty} H_n \left( \frac{1}{n} - \frac{1}{n+1} \right) = \sum_{n=1}^{\infty} \left( \frac{H_n}{n} - \frac{H_{n+1}}{n+1} + \frac{1}{(n+1)^2} \right) = \varsigma(2)$$

Hence we obtain

$$\int_0^1 \xi(x)\, dx = \sum_{n=1}^{\infty} H_n \left( \log\frac{n+1}{n} - \frac{1}{n} \right) + \varsigma(2)$$

and using (9.1) this becomes

$$2\int_0^1 \xi(x)\, dx = \varsigma(2) - \gamma^2 + 2\gamma_1$$

From (9.4) we have

$$2\int_0^1 \xi(x)\, dx = \int_0^1 \left( \psi'(x) - [\psi(x)+\gamma]^2 \right) dx - \varsigma(2)$$

$$\int_0^1 \left( \psi'(x) - [\psi(x)+\gamma]^2 \right) dx - \varsigma(2) = \varsigma(2) - \gamma^2 + 2\gamma_1$$

and this gives us an apparently new integral representation for $\gamma_1$

$$(9.5) \qquad \int_0^1 \left( \psi'(x) - \psi^2(x) - 2\gamma\psi(x) \right) dx = 2\varsigma(2) + 2\gamma_1$$

Incidentally, differentiating Nielsen's formula (9.4) gives us

$$(9.6) \qquad 2[\psi(x) + \gamma]\psi'(x) = \psi''(x) + 2\sum_{n=1}^{\infty} \frac{H_n}{(x+n)^2}$$

and letting $x = 1$ results in

$$\psi''(1) = 2\sum_{n=1}^{\infty} \frac{H_n}{(n+1)^2} = 2\sum_{n=1}^{\infty} \frac{H_{n+1}}{(n+1)^2} - 2\sum_{n=1}^{\infty} \frac{1}{(n+1)^3}$$

$$= 2\sum_{n=1}^{\infty} \frac{H_n}{n^2} - 2\varsigma(3)$$

Since [sr, p.22]

$$\psi^{(n)}(x) = (-1)^{n+1} n! \varsigma(n+1, x)$$

this gives us the well-known Euler sum

$$\sum_{n=1}^{\infty} \frac{H_n}{n^2} = 2\varsigma(3)$$

Using $H_n = \psi(n+1) + \gamma$ we may write (9.6) as

$$[\psi(x) + \gamma]\varsigma(2, x) + \varsigma(3, x) = \sum_{n=1}^{\infty} \frac{\psi(n+1)}{(x+n)^2} + \gamma \sum_{n=1}^{\infty} \frac{1}{(x+n)^2}$$

$$= \sum_{n=1}^{\infty} \frac{\psi(n+1)}{(x+n)^2} + \gamma\varsigma(2, x) - \frac{\gamma}{x^2}$$

This gives us

$$\psi(x)\varsigma(2, x) + \varsigma(3, x) = \sum_{n=1}^{\infty} \frac{\psi(n+1)}{(x+n)^2} - \frac{\gamma}{x^2}$$

or equivalently

$$\psi(x)\varsigma(2,x) + \varsigma(3,x) = \sum_{n=0}^{\infty} \frac{\psi(n+1)}{(x+n)^2} \tag{9.7}$$

With $x = 1$ we have

$$-\gamma\varsigma(2) + \varsigma(3) = \sum_{n=0}^{\infty} \frac{\psi(n+1)}{(n+1)^2} = \sum_{n=1}^{\infty} \frac{\psi(n)}{n^2}$$

as previously noted by Coffey [18] (and see (8.6) above).

It is clear that $\sum_{n=0}^{\infty} \frac{\psi(n+1)}{(x+n)^p}$ may be evaluated by differentiating (9.6) $p - 2$ times.

Donal F. Connon
Elmhurst
Dundle Road
Matfield
Kent TN12 7HD
dconnon@btopenworld.com